%% file: main.tex
\documentclass[10pt,journal]{IEEEtran}
\IEEEoverridecommandlockouts
\usepackage{cite}
\usepackage{amsmath,amssymb,amsfonts}
\usepackage{graphicx}
\usepackage{amsmath}
\usepackage{xcolor}
\usepackage{colortbl}
\renewcommand{\angle}{\measuredangle}
\newcommand\myworries[1]{\textcolor{blue}{#1}}
\usepackage{amsthm}
\newtheorem{remark}{Remark} %\newtheorem* --> for no number of remark.

\newtheorem*{theorem*}{Theorem}

\newtheorem*{lemma*}{Lemma}

\usepackage{steinmetz}
\newtheorem*{definition*}{Definition}

\newtheorem*{condition*}{Condition}
\newtheorem*{proof*}{Proof}

\usepackage{algorithm}% http://ctan.org/pkg/algorithm
\usepackage{algpseudocode}% http://ctan.org/pkg/algorithmicx

\makeatletter
\def\therule{\makebox[\algorithmicindent][l]{\hspace*{.5em}\vrule height .75\baselineskip depth .25\baselineskip}}%

\newtoks\therules% Contains rules
\therules={}% Start with empty token list
\def\appendto#1#2{\expandafter#1\expandafter{\the#1#2}}% Append to token list
\def\gobblefirst#1{% Remove (first) from token list
	#1\expandafter\expandafter\expandafter{\expandafter\@gobble\the#1}}%
\def\LState{\State\unskip\the\therules}% New line-state
\def\pushindent{\appendto\therules\therule}%
\def\popindent{\gobblefirst\therules}%
\def\printindent{\unskip\the\therules}%
\def\printandpush{\printindent\pushindent}%
\def\popandprint{\popindent\printindent}%

\algdef{SE}[WHILE]{While}{EndWhile}[1]
{\printandpush\algorithmicwhile\ #1\ \algorithmicdo}
{\popandprint\algorithmicend\ \algorithmicwhile}%
\algdef{SE}[FOR]{For}{EndFor}[1]
{\printandpush\algorithmicfor\ #1\ \algorithmicdo}
{\popandprint\algorithmicend\ \algorithmicfor}%
\algdef{S}[FOR]{ForAll}[1]
{\printindent\algorithmicforall\ #1\ \algorithmicdo}%
\algdef{SE}[LOOP]{Loop}{EndLoop}
{\printandpush\algorithmicloop}
{\popandprint\algorithmicend\ \algorithmicloop}%
\algdef{SE}[REPEAT]{Repeat}{Until}
{\printandpush\algorithmicrepeat}[1]
{\popandprint\algorithmicuntil\ #1}%
\algdef{SE}[IF]{If}{EndIf}[1]
{\printandpush\algorithmicif\ #1\ \algorithmicthen}
{\popandprint\algorithmicend\ \algorithmicif}%
\algdef{C}[IF]{IF}{ElsIf}[1]
{\popandprint\pushindent\algorithmicelse\ \algorithmicif\ #1\ \algorithmicthen}%
\algdef{Ce}[ELSE]{IF}{Else}{EndIf}
{\popandprint\pushindent\algorithmicelse}%
\algdef{SE}[PROCEDURE]{Procedure}{EndProcedure}[2]
{\printandpush\algorithmicprocedure\ \textproc{#1}\ifthenelse{\equal{#2}{}}{}{(#2)}}%
{\popandprint\algorithmicend\ \algorithmicprocedure}%
\algdef{SE}[FUNCTION]{Function}{EndFunction}[2]
{\printandpush\algorithmicfunction\ \textproc{#1}\ifthenelse{\equal{#2}{}}{}{(#2)}}%
{\popandprint\algorithmicend\ \algorithmicfunction}%
\makeatother

\renewcommand\myworries[1]{}
\usepackage{mathtools}

\makeatletter
\newsavebox\myboxA
\newsavebox\myboxB
\newlength\mylenA

\newcommand*\xoverline[2][0.75]{%
	\sbox{\myboxA}{$\m@th#2$}%
	\setbox\myboxB\null% Phantom box
	\ht\myboxB=\ht\myboxA%
	\dp\myboxB=\dp\myboxA%
	\wd\myboxB=#1\wd\myboxA% Scale phantom
	\sbox\myboxB{$\m@th\overline{\copy\myboxB}$}%  Overlined phantom
	\setlength\mylenA{\the\wd\myboxA}%   calc width diff
	\addtolength\mylenA{-\the\wd\myboxB}%
	\ifdim\wd\myboxB<\wd\myboxA%
	\rlap{\hskip 0.5\mylenA\usebox\myboxB}{\usebox\myboxA}%
	\else
	\hskip -0.5\mylenA\rlap{\usebox\myboxA}{\hskip 0.5\mylenA\usebox\myboxB}%
	\fi}
\makeatother

\usepackage{scalerel}[2014/03/10]
\usepackage{stackengine}

\newlength\myindent
\usepackage{epsfig}
\usepackage{tabularx}
\usepackage{supertabular,booktabs}
\usepackage{capt-of}

\usepackage{longtable}
\usepackage{siunitx}
\usepackage{booktabs}

\usepackage{subcaption}

\usepackage{algorithm}
\usepackage{algpseudocode}
\usepackage{bm}

\renewcommand{\algorithmicforall}{\textbf{for each}}

\newcommand{\bd}{\mathbf}

\algrenewcommand\algorithmicrequire{\textbf{Precondition:}}
\algrenewcommand\algorithmicensure{\textbf{Postcondition:}}

\usepackage{textcomp}
\usepackage[normalem]{ulem}
\usepackage{mathtools}
\usepackage{float}
\usepackage{tikz}
\usetikzlibrary{shapes.geometric, arrows}
\tikzstyle{startstop} = [rectangle, rounded corners, minimum width=3cm, minimum height=1cm,text centered, draw=black, fill=red!30]
\tikzstyle{io} = [trapezium, trapezium left angle=70, trapezium right angle=110, minimum width=3cm, minimum height=1cm, text centered, text width=3cm, draw=black, fill=blue!30]
\tikzstyle{process} = [rectangle, minimum width=3cm, minimum height=1cm, text centered, draw=black, fill=orange!30]
\tikzstyle{decision} = [diamond, minimum width=2.5cm, minimum height=2.5cm, text centered, draw=black, fill=green!30]
\tikzstyle{arrow} = [thick,->,>=stealth]
\tikzstyle{process} = [rectangle, minimum width=3cm, minimum height=1cm, text centered, text width=3cm, draw=black, fill=orange!30]
\def\BibTeX{{\rm B\kern-.05em{\sc i\kern-.025em b}\kern-.08em
		T\kern-.1667em\lower.7ex\hbox{E}\kern-.125emX}}
\begin{document}

	\title{Power Flow as Intersection of Circles:\\ A New Fixed Point Method \thanks{This work is partially supported by the National Science Foundation awards ECCS-1810537 and ECCS-1807142.} \thanks{K. P. Guddanti and Y. Weng are with the School of Electrical, Computer and Energy Engineering at Arizona State University, emails:\{kguddant,yang.weng\}@asu.edu; B. Zhang is with the Department of Electrical and Computer Engineering at the University of Washington, email: zhangbao@uw.edu.}
	}
	\author{
		\IEEEauthorblockN{Kishan Prudhvi Guddanti},
		\IEEEauthorblockA{\textit{Student Member}},
		\textit{IEEE},
		\and
		\IEEEauthorblockN{Yang Weng},
		\IEEEauthorblockA{\textit{Member}},
		\textit{IEEE},
		\and
		\IEEEauthorblockN{Baosen Zhang},
		\IEEEauthorblockA{\textit{Member}},
		\textit{IEEE}
	\vspace{-10mm}}
	\DeclarePairedDelimiter\abs{\lvert}{\rvert}
	\makeatletter
	\let\oldabs\abs
	\def\abs{\@ifstar{\oldabs}{\oldabs*}}
	\maketitle

	\begin{abstract}
The power flow (PF) problem is a fundamental problem in power system engineering. Many popular solvers face challenges, such as convergence issues. One can try to rewrite the PF problem into a fixed point equation, which can be solved exponentially fast. But, existing methods have their own restrictions, such as the required AC network structure or bus types. To remove these restrictions, we employ the circle geometry per-bus via rectangular coordinate representation to embed our physical knowledge of operation point selection in PV curves. Each iteration of the algorithm consists of finding intersections of circles, which can be computed efficiently with high numerical accuracy. Such analysis also helps in visualizing PV curve to always select the high voltage solution. We compare the performance of our fixed point algorithm with existing state-of-the-art methods, showing that the proposed method can correctly find the solutions when other methods cannot. In addition, we empirically show that the fixed point algorithm is much more robust to bad initialization points than the existing methods.
%leading to simply repeatition method to solve the problems.   were proposed with a unique property: If   
%For such problems, we present a new class of algorithms obtained from writing the power flow problem as finding the solution of a fixed point equation. These fixed point equations can be solved iteratively by simple function evaluations and can potentially converge to the solution exponentially fast. 
%Our approach removed such limits by a new geometric insight into the power flow equations under rectangular coordinates. 
% Different than existing fixed point based algorithms, our proposed method encompasses  The Jacobian-based solvers are widely accepted and employed due to their advantage in speed of solving the PF problem in very few iterations. However, these solvers may become unstable for ill-conditioned systems. This is because of the numerical issue of highly sensitive Jacobian when nearing singularity. This paper introduces a new fixed-point representation of power flow equations that can be used to solve such ill-conditioned systems. The fixed-point equation is derived using the representation of power flow equations as circles. The proposed method can include all possible cases: PV/PQ buses, mesh network, fully resistive and inductive lines. It is also observed that the proposed rectangular fixed-point power flow (RFPPF) solver has exponential convergence characteristics. The proposed method is also insensitive to random initialization in the case of a network with PQ buses.
\end{abstract}


	\section{Introduction}
	\input{intro.tex}
	\section{Power Flow Equations and Circles}
	\label{sec:pf}
    \input{powerflow.tex}

	\section{Fixed Point Equation for Power Flow}
	\label{sec:fixed_point}
	 \input{fp.tex}

	 \section{Main Algorithm}
	 \label{sec:algorithm}
	 \input{algorithm.tex}

	\subsection{Finding Intersection Points}
	\label{sec:FPE_powerflow_prob}
	\input{fixed_point_equation.tex}
	\section{Numerical Results}
	\label{sec:results}
	\input{numerical.tex}

	\section{Conclusion}
	\label{sec:conclu}
	\input{conclusion.tex}

    \bibliographystyle{IEEEtran}
    \bibliography{Kishan,bz}
\vspace{-1em}
    \appendix
    \label{sec:app}
    \input{app.tex}
\begin{IEEEbiography}[{\includegraphics[width=1in,height=1.25in,clip,keepaspectratio]{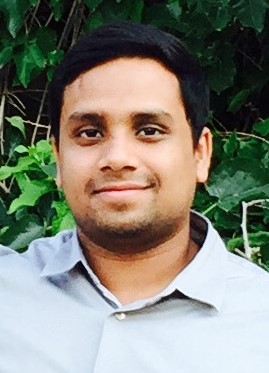}}]{Kishan Prudhvi Guddanti} received the B.Tech. degree in electrical engineering from Sri Ramaswamy Memorial Institute of Science and Technology, Chennai, India and the M.Sc degree in electrical engineering from Arizona State University, Tempe, AZ, USA.

		He is currently pursuing the Ph.D. degree at Arizona State University, Tempe, AZ, USA. He is one of the winners of an AI competition organized by RTE France in 2019. His current research interest is in the interdisciplinary area of AI applications in power systems in addition to voltage stability and data-driven techniques applications in power system risk assessment and control.
\end{IEEEbiography}

\begin{IEEEbiography}[{\includegraphics[width=1in,height=1.25in,clip,keepaspectratio]{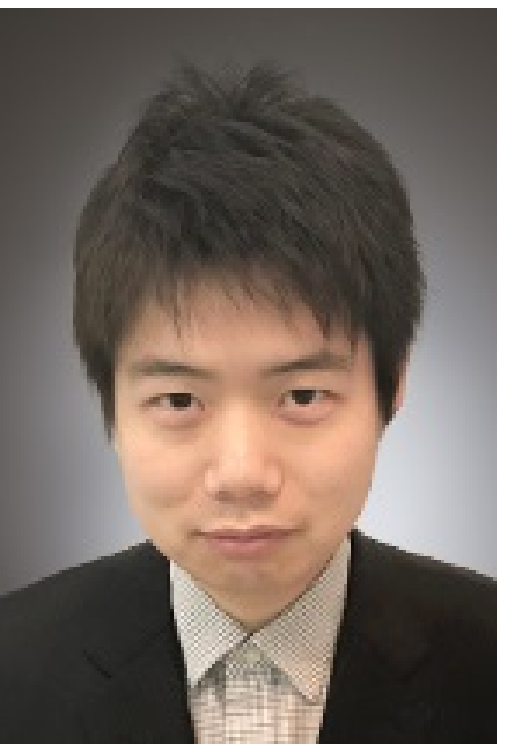}}]{Yang Weng}
received the B.E. degree in electrical engineering from Huazhong University of Science and Technology, Wuhan, China; the M.Sc. degree in statistics from the University of Illinois at Chicago, Chicago, IL, USA; and the M.Sc. degree in machine learning of computer science and M.E. and Ph.D. degrees in electrical and computer engineering from Carnegie Mellon University (CMU), Pittsburgh, PA, USA.

After finishing his Ph.D., he joined Stanford University, Stanford, CA, USA, as the TomKat Fellow for Sustainable Energy. He is currently an Assistant Professor of electrical, computer and energy engineering at Arizona State University (ASU), Tempe, AZ, USA. His research interest is in the interdisciplinary area of power systems, machine learning, and renewable integration.

Dr. Weng received the CMU Dean's Graduate Fellowship in 2010, the Best Paper Award at the International Conference on Smart Grid Communication (SGC) in 2012, the first ranking paper of SGC in 2013, Best Papers at the Power and Energy Society General Meeting in 2014, ABB fellowship in 2014, Golden Best Paper Award at the International Conference on Probabilistic Methods Applied to Power Systems in 2016, and Best Paper Award at IEEE Conference on Energy Internet and Energy system Integration in 2017
\end{IEEEbiography}
\vspace{-38em}
\begin{IEEEbiography}
	[
	{\includegraphics[width=1in,height=1.25in,clip,keepaspectratio]{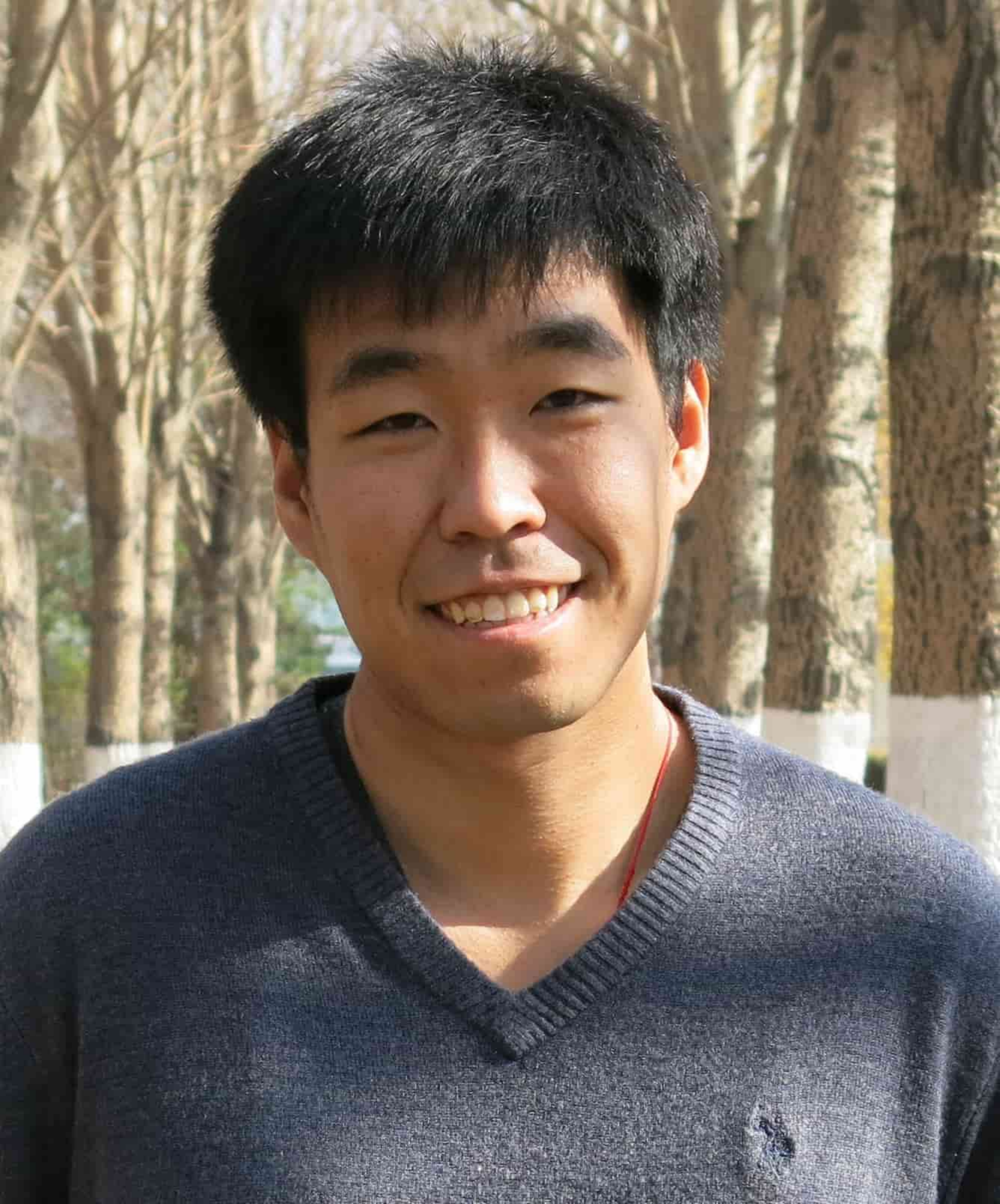}}
	] {Baosen Zhang} received his Bachelor of Applied Science in Engineering Science degree from the University of Toronto in 2008; and his PhD degree in Electrical Engineering and Computer Sciences from University of California, Berkeley in 2013.

	He was a Postdoctoral Scholar at Stanford University, affiliated with the Civil and Environmental Engineering and Management \& Science Engineering. He is currently the Keith and Nancy Rattie Endowed Career Development Professor in the Departement of Electrical and Computer Engineering at the University of Washington, Seattle, WA. His research interests are in power systems and cyberphysical systems.
\end{IEEEbiography}

\end{document}

%% file: intro.tex
%% !TEX root=main.tex
The power flow problem is one of the canonical problems in power engineering and it is frequently used in power system operation and planning studies~\cite{Stott74,Glover}.
Existing power flow methods mostly rely on iterative methods such as Newton-Raphson (NR)~\cite{Momoh99} or fast decoupled load flow (FDLF)~\cite{stott1974fast,MonticelliEtAl1990}. These algorithms have been the workhorses of the power industry and have performed well most of the time.  
% f Despite being the workhorses of power system analysis,
However, as large-scale development of renewable resources and distributed generation push systems to operate in new regimes, the existing algorithms can experience convergence issues, especially when systems operate close to their loadability limits~\cite{Thorp89,Thorp90,Rao2018}. 
% these algorithms can suffer from convergence issues, especially when the system operates close to its loadability limits or the initial guess is poor~\cite{Thorp89,Thorp90}. Both of these issues are exacerbated by the large-scale development of renewable resources and distributed generation and 
Therefore, the need for new efficient and robust power flow algorithms to complement these existing methods remains despite decades of studies~\cite{pudjianto2007virtual}.

% The power flow equations are relaxed to solve the power flow problem faster with an approximate solution using DC power flow solver \cite{Stott74}. However, the DC power flow solver is not reliable in every case because of its relaxed form. For example, in case of distribution system where the resistance to reactance ratio ($R/X$) is high, the DC power flow solver doesn't provide a sufficiently accurate solution \cite{Rajicic88}, it also has a poor convergence \cite{Zimmerman95}. It is important to have a power flow solver that has reliable performance on both transmission and distribution side \cite{Singhal17}. On the other hand, the AC power flow solver can solve the power flow problem on both the transmission and distribution side effectively. In order to solve the AC power flow problem faster, the Jacobian-based power flow solver is used.

Algorithms like NR can be thought as variants of descent algorithms (or approximate descent in the case of FDLF) that modifies the solution iteratively. A fundamental reason for why these algorithms can fail to converge to a solution is simply because the geometry of power flow is not convex~\cite{Zhang13,Hiskens01}.  For example, NR uses the Jacobian to find the direction of the steepest descent. Because the power flow equations are nonlinear and nonconvex, there are many local minimums and saddle points, and the Jacobian fails to provide a meaningful descent direction at these points. To prevent the algorithm from getting stuck, it becomes important to pick "good" initial starting points~\cite{Ajjarapu92, Stot71, Iwa78, Canizares93}. Consequently, a number of methods have been developed to overcome the sensitive dependence on the initial guess \cite{Chiang14, Milano09, Costa08,Exposito13}.
% \textcolor{blue}{Specifically, Jacobian incorporates information of all buses into a single matrix and provides the direction of voltage update for all buses at every iteration of NR method. Even though this voltage update solves the problem faster, it may or may not lead to the solution when the Jacobian is closer to singularity or a ``bad'' initial guess is provided to the solver. This is because the Jacobian fails to take a meaningful voltage update towards the solution.} Recognizing this, a series of studies have used techniques such as homotopy, current injection and factorized methods to overcome the sensitive dependence on the initial guess \cite{Chiang14, Milano09, Costa08,Exposito13}.

As systems start to operate closer to there limits, picking better initialization points becomes insufficient. Since the Jacobians for all points close to the boundary of the feasible power flow region have eigenvalues close to 0 (they loose rank), they necessarily become ill-conditioned and iterative algorithms may diverge~\cite{Tinney67,Shah12,Wood12}.  
% are all poorly conditioned have Ja A second challenge in adopting NR and methods comes from the nature of the power flow Jacobian, where Jacobian provides update directions in iterations. However, when the Jacobian becomes singular or ill-conditioned, e.g., when the loading is heavy~\cite{Tinney67,Shah12,Wood12}, the Jacobian-based algorithms tend to diverge. 
To avoid this phenomenon, a class of non-divergent power flow algorithms was developed to accelerate or decelerate the updates based on the conditioning of the Jacobian~\cite{Iwamoto81,Castro97,Braz2000,Bijwe03}. However, these approaches can still be sensitive to the initial guess and sometimes exhibit oscillatory behavior, where the solutions may neither converge nor diverge. An approach using complementarity conditions is developed in~\cite{Pirnia13}, but it can reach local minimums or saddle points instead of the true power flow solution. Energy-bases analysis based on mechanical models can help algorithms to escape these stationary points~\cite{moon2003energy}, but implementing them for different bus types in a practical power system is non-trivial. 
% \textcolor{blue}{To avoid this problem i.e., inversion of ill-conditioned Jacobian, optimization framework is used to solve the power flow problem~\cite{Pirnia13}. However, these methods suffer from getting stuck at saddle points or local optimums due to non-linear nature of the power flow problem~\cite{Pirnia13}. \cite{moon2003energy} overcomes the problem of non-linear power flow optimization by using energy-based analysis that use mechanical model of power system. However, the application of bus types and other devices into these energy functions is nontrivial.}

Recently, a new class of power flow formulations based on fixed point equations has been proposed to overcome the algorithmic challenges present in descent algorithms~\cite{John17,Wang18}. The basic idea is to write the power flow equations in a form of $\bd v= f(\bd v)$, where $\bd v$ is the complex voltage and a fixed point of the function $f$~\cite{kakutani1941generalization}. If this relationship can be found, then a simple algorithm to find the fixed point is to repeatedly apply the function $f$. Furthermore, if the iterates converge to the fixed point, then it will converge exponentially quickly. The challenge is to find a suitable $f$, which have only exists for restricted class of systems. For example, the results in~\cite{Wang18} applies to networks with only PQ buses and the results in~\cite{John17} only applies to purely inductive (lossless) radial networks.

% \textcolor{blue}{However, it would be easy to take a meaningful voltage update step even when the system is ill-conditioned by visualizing the PV curve. By observing the PV curve at every bus, we can always say that we want the high voltage solution when the solver is solving the problem. Hence a per-bus analysis can help to solve this problem as it allows to visualize the PV curve at every bus and select the high voltage solution at every bus. For an $n$-bus system, usually network reduction techniques are employed to visualize the PV curve \cite{Vu99} but for a power flow problem, we cannot afford to do network reduction. Hence in this paper, a new mathematical formulation is provided for the first time to help visualize the PV curve for $n$-bus power system without using any network reduction techniques and always select the high voltage solution.}

%we show that the problem of meaningful voltage update for an ill-conditioned system can be understood easily by visualizing the per-bus PV curve in our proposed model. By observing the PV curve at every bus, we can actively select the high voltage solution when solving the problem. 

\begin{figure*}[t!]
		\centering
        		\begin{subfigure}{0.34\textwidth}
			%\vspace{3.5mm}
		\includegraphics[scale=0.5]{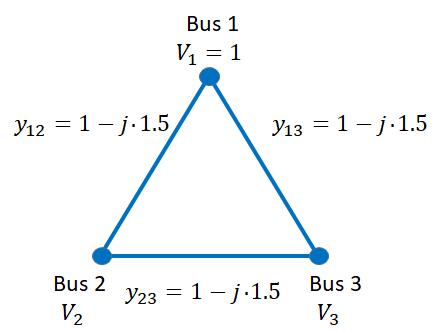}
		\caption{3-bus system.}
		\label{3busfig}
		\end{subfigure}
        ~
		\begin{subfigure}{0.315\textwidth}
			\vspace{1.5mm}
			\includegraphics[width=0.7\textwidth]{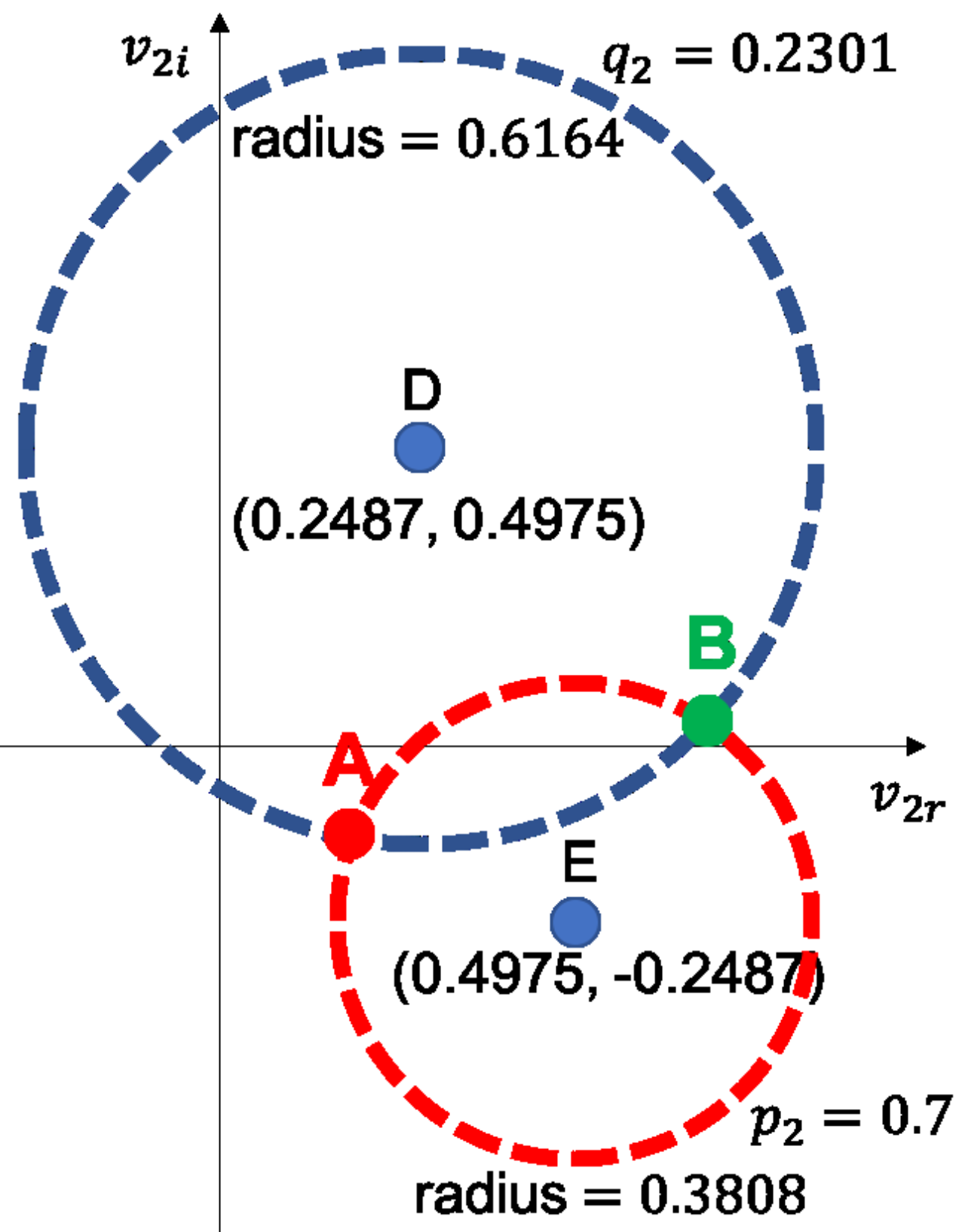}
			%\newline
			\vspace{0mm}
			\caption{P and Q at bus 2.}
			\label{fig:bus2}
		\end{subfigure}
		~
		\begin{subfigure}{0.22\textwidth}
			\includegraphics[width=\textwidth]{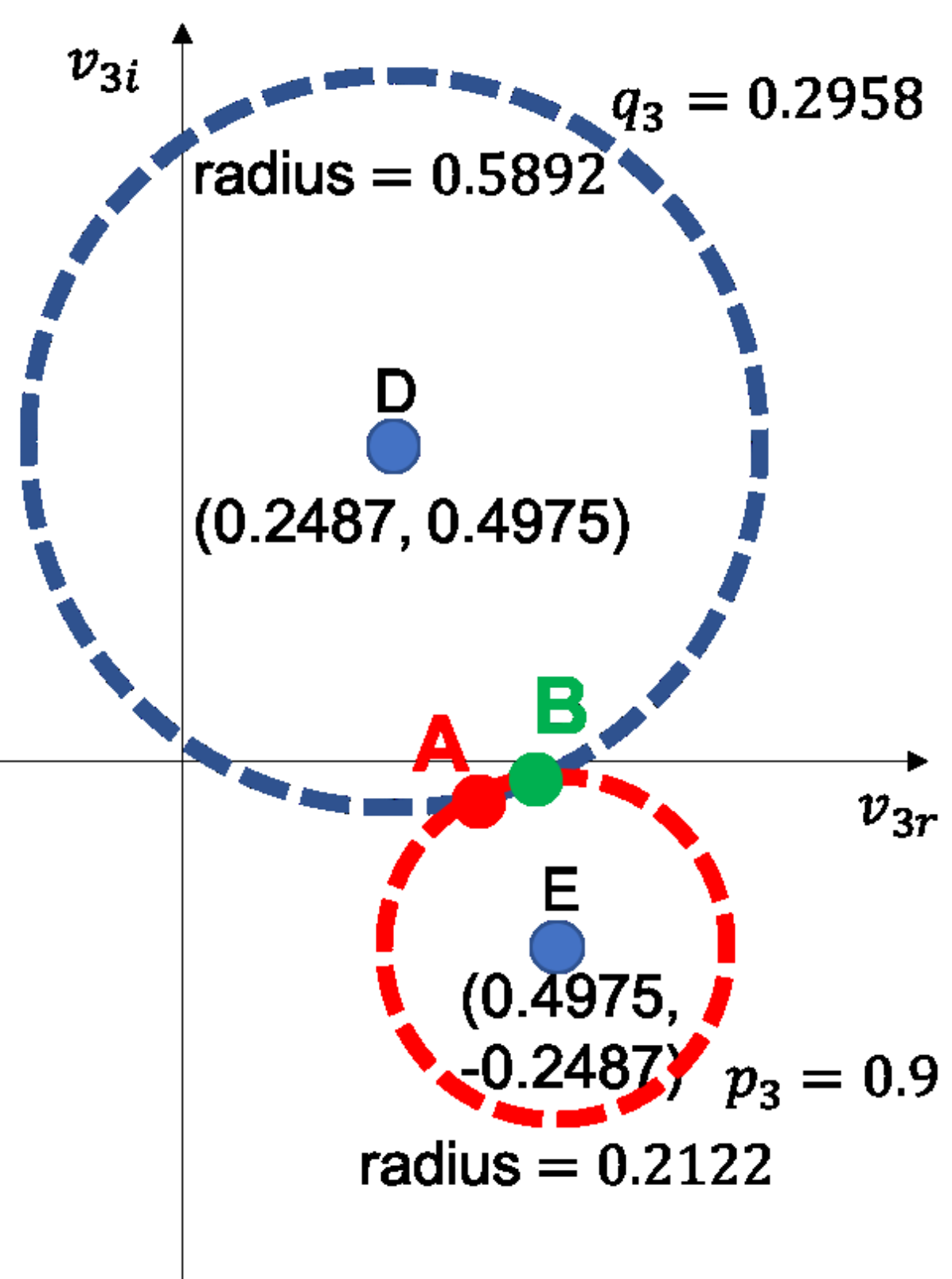}
			\vspace{-6mm}
			\caption{P and Q at bus 3.}
			\label{fig:bus3}
		\end{subfigure}
		% \vspace{-2mm}
		\caption{Active and reactive circles for a three-bus system. Bus 1 is the slack bus and buses 2 and 3 are PQ buses.}
		\label{fig:fixedpoint}
		%\vspace{-7mm}
 \end{figure*}
 
In this paper, we present a novel fixed point formulation of the full AC power flow equations that is applicable to networks with \emph{arbitrary topologies and mixture of PQ and PV buses}.   
% propose to solve power flow problem based on a fixed point method, where the fixed point is defined as a vector of complex voltages. Such 
This approach work is based on a coordinate transformation, where power flow solutions are interpreted as the intersections of circles, where the parameters (center and radius) of the circles depend linearly on the voltages of the neighboring buses. 
% derived and generalized for each bus via mathematically intersected per-bus circles. % are consistent with the parameters of the neighbors. 
This formulation can be thought as a generalization of the PV noise curve often used to visualize power transfer between two buses. Computationally, only the intersection of two circles needs to be calculated, which involves a series of simple algebraic computations. Therefore, this approach is much cheaper than algorithms~(e.g., NR) that require matrix calculations.

% provided for the first time to help visualize and utilize the per-bus PV curve for $n$-bus power flow analysis. For example, our fixed point formulation view the power flow solution as the \emph{intersection of circles}, where the parameters (center and radius) of the circles depend linearly on the voltages of the neighboring buses.} During each iterate, only the intersection of two circles needs to be calculated, which involves a series of simple algebraic computations. Therefore, this step is much cheaper than algorithms~(e.g., NR) that require matrix calculations. It is important to highlight that %even though we suspect that our method is a contraction mapping, due to the somewhat complicated dependence between the buses, we cannot provide a theoretical guarantee of its convergence. This is in contrast to the theoretical results in~\cite{John17,Wang18}, which have less generality and stricter assumptions, but do provide some guarantees on convergence. 
% \emph{our formulation is valid for full AC networks with arbitrary topology and arbitrary mix of PV and PQ buses, previously not exact in the literature}. This comes without using any network reduction techniques and provides the flexibility to select the high voltage solution during iterations.

To verify the performance of our algorithm, we test it on the standard IEEE systems, including large ones with $2383$ and 3375 buses. We compare our approach with NR, FDLF and non-divergent power flow algorithms. We show that when the loading is heavy, our algorithm is able to converge to the right solution while the other algorithm can diverge or become unstable. In addition, we show that our method is much more robust to random initialization points than the other methods. It is important to note that we are not advocating to replace existing power flow solvers. These algorithms have been highly optimized and do perform extremely well in many situations. Rather, the proposed algorithm in this paper can be used as a complementary tool by the system operators when conventional algorithm diverge or stall.  

The paper is organized as follows: Section~\ref{sec:pf} introduces the rectangular power flow equations and show how they can be thought as intersections of circles. Section~\ref{sec:fixed_point} discusses the fixed point formulation of the power flow equations and walks through a three-bus example. Section~\ref{sec:algorithm} presents the main algorithm. Section~\ref{sec:FPE_powerflow_prob} introduces a $3$-tuple vector form of circles and shows how closed-form formulas with good numerical properties can be found using the vector notation. Section~\ref{sec:results} shows numerical results of our proposed algorithm compared against existing state-of-the-art algorithms on different IEEE benchmark networks. Section~\ref{sec:conclu} concludes the paper.

%% file: powerflow.tex
%% !TEX root=main.tex

\subsection{Power Flow Equations in Rectangular Coordinates}
	\label{intro_rec_pf_eqs}
	%The geometric understanding of rectangular coordinate-based  power flow equations for loadability analysis is discussed in \cite{Weng}. These rectangular coordinate-based power flow equations have three main advantages. $1)$ Reduction of complexity by converting nonlinear equations with polynomials and sinusoids to purely polynomial. 2) Power flow is solved in this paper with linear updates which avoids the singularity of Jacobian. 3) The proposed algorithm is faster when compared with the traditional Newton-Raphson method.

Throughout this paper we use rectangular coordinates where a bus is index by $d$; $p_d$ and $q_d$ are the active and reactive powers, respectively; $\mathit{v}_{d,r}$ and $\mathit{v}_{d,i}$ are the real and imaginary parts of the bus voltage, respectively; and $\mathcal{N}(d)$ is the set of neighboring buses connected to bus $d$. We adopt the standard $\Pi$ model of transmission lines~\cite{Glover} and write the admittance of a line between buses $d$ and $k$ as $g_{dk}+jb_{dk}$. We assume that $b_{dk} \leq 0$ for all lines~(lines are inductive). In these notations, the power flow equations become~\cite{TorresEtAl1998,weng2017geometric}:
\begin{equation}
	\mathit{p}_d
	= t_{d,1} \cdot \mathit{v}^2_{d,r} + t_{d,2} \cdot \mathit{v}_{d,r} + t_{d,1} \cdot \mathit{v}^2_{d,i} + t_{d,3} \cdot \mathit{v}_{d,i} ,
	\label{math_p}
	\end{equation}
	\begin{equation}
	\mathit{q}_d
	= t_{d,4} \cdot \mathit{v}^2_{d,r } - t_{d,3} \cdot \mathit{v}_{d,r} + t_{d,4} \cdot \mathit{v}^2_{d,i} + t_{d,2} \cdot \mathit{v}_{d,i}.
	\label{math_q}
	\end{equation}
The parameters $t_{d,1},t_{d,2},t_{d,3},t_{d,4}$ are given by
	\begin{equation*}
	t_{d,1}
	= - \sum_{\mathit{k}\in\mathcal{N}(\mathit{d})} \mathit{g}_{k,d} ,
	\ t_{d,2}
	= \sum_{\mathit{k}\in\mathcal{N}(\mathit{d})} (\mathit{v}_{k,r}\mathit{g}_{k,d} - \mathit{v}_{k,i}\mathit{b}_{k,d}) ,
	\label{ref1}
	\end{equation*}
	\begin{equation*}
	t_{d,3}
	= \sum_{\mathit{k}\in\mathcal{N}(\mathit{d})} (\mathit{v}_{k,r}\mathit{b}_{k,d} + \mathit{v}_{k,i}\mathit{g}_{k,d}) ,
	\ t_{d,4}
	=  \sum_{\mathit{k}\in\mathcal{N}(\mathit{d})} \mathit{b}_{k,d}.
	\label{ref2}
	\end{equation*}
The validity of the above equations can be checked by straightforward substitutions.

Since the terms $t_{d,1}$ and $t_{d,4}$ are always negative, \eqref{math_p} and \eqref{math_q} describe two circles in the variables $v_{d,r}$ and $v_{d,i}$. We call the circle described by \eqref{math_p} the active power circle parametrized by its center $\mathbb{C}_{p}$ and radius $\mathit{r}_{p}$; similarly, we say that \eqref{math_q} describes the reactive power circle parameterized by center $\mathbb{C}_{q}$ and radius $\mathit{r}_{q}$. These parameters are given by:
\begin{align}
&\mathbb{C}_{p}
= \left(\frac{-t_{d,2}}{2t_{d,1}}, \frac{-t_{d,3}}{2t_{d,1}}\right), \;
\mathbb{C}_{q}
= \left(\frac{t_{d,3}}{2t_{d,4}}, \frac{-t_{d,2}}{2t_{d,4}}\right), \label{eq:cpq}\\
&\mathit{r}_{p}
= \sqrt{\frac{\mathit{p}_{d}}{t_{d,1}} +  \frac{\left(t_{d,2}\right)^2 + \left(t_{d,3}\right)^2}{4t^2_{d,1}}}, \label{eq:rp}\\
&\mathit{r}_{q}
= \sqrt{\frac{\mathit{q}_d}{t_{d,4}} +  \frac{\left(t_{d,3}\right)^2 + \left(t_{d,2}\right)^2}{4t^2_{d,4}}}.\label{eq:rq}
\end{align}

Figure~\ref{fig:fixedpoint} shows a three bus network and the associated active and reactive power cycles are buses 2 and 3~(bus 1 is assumed to the slack bus). The intersection points A and B in Fig.~\ref{fig:bus2} and Fig.~\ref{fig:bus3} represents the potential power flow solutions.
% \begin{figure*}[ht!]
% 		\centering
%     \begin{subfigure}{0.27\textwidth}
% 			%\vspace{3.5mm}
% 		\includegraphics[scale=1]{full_graph.png}
% 		\caption{6 bus system.}
% 		\label{fig:6bus}
% 		\end{subfigure}
%         ~
% 		\begin{subfigure}{0.7\textwidth}
% 			\vspace{1.5mm}
% 			\includegraphics[width=\textwidth,height=4cm]{update_graph_rule.png}
% 			%\newline
% 			\vspace{0mm}
% 			\caption{Illustration of the update process in one round.}
% 			\label{fig:update_order}
% 		\end{subfigure}
% 		\caption{An example of a six bus system where the buses are updated in lexicographical order in each round.}
% 		\label{fig:update}
% 		%\vspace{-7mm}
%  \end{figure*}
 
    \begin{figure*}[ht!]
		\centerline{\includegraphics[width=\linewidth]{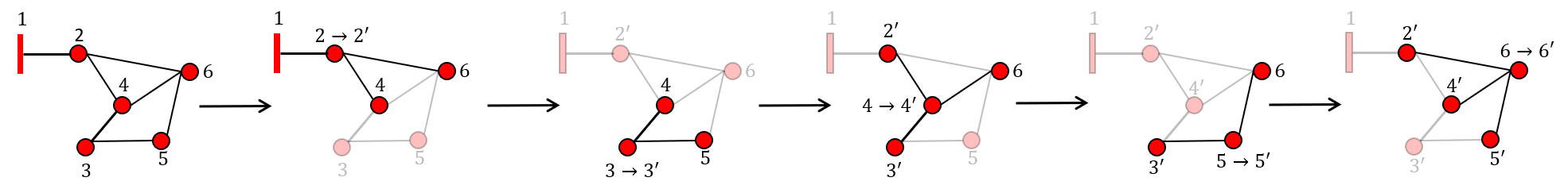}}
		\caption{An example of a six bus system where the buses are updated in lexicographical order in one round.}
		\label{fig:update}
    \end{figure*}
 
\subsection{PV Buses} \label{subsec:PV_bus}
The discussions in the above section focus on PQ buses, but PV buses are also frequently used to describe generators~\cite{GarciaEtAl2000}. In this case, the reactive power balance equation in~\eqref{math_q} is replaced by a condition on the voltage magnitude:
\begin{equation} \label{eqn:voltage}
v_{d,r}^2+v_{d,i}^2 = V_{ref}^2,
\end{equation}
where $V_{ref}$ is the reference voltage. Again, we can think of PV buses in term of circles, since \eqref{eqn:voltage} is a circle centered at the origin with a fixed radius. Therefore, our framework does not require different treatment of PQ and PV buses.

%% file: fp.tex
%% !TEX root=main.tex

The geometric representation of the power flow equations as the intersection of circles leads to a simple fixed point view of power flow solutions. Suppose that a vector of complex voltages is given. Then, the voltage at a particular bus $d$ is determined by its neighbors as the intersection of the active power circle with the reactive power circle (for a PQ bus) or with the voltage magnitude circle (for a PV bus). Of course, two circles, if they intersect, could do so at two distinct points as shown in Figs.~\ref{fig:bus2} and~\ref{fig:bus3}. In this case, we need to pick one of the intersection points as the complex voltage at a bus and use it to compute the parameter of its neighboring circles. 
% A natural question arises: if the circles intersect in two points as shown in Figs.~\ref{fig:bus2} and~\ref{fig:bus3}, which of the point should we choose as the voltage at that bus? 
To make this choice, we follow two common assumptions made in power flow calculations.

The first assumption we make is that we are interested in solutions at higher voltage magnitudes~\cite{DOBSON1989253,kwatny1995local,Johnson77}. These solutions have long been seen as the practical and stable solutions in actual systems~\cite{Mehta16}. For example, in both Figs.~\ref{fig:bus2} and~\ref{fig:bus3}, we would chose point $B$ as the solution. For a PV bus, all points of intersection have the same voltage magnitude. In this case, we make the second assumption that voltages with smaller (absolute) angles are preferable. This assumption is rooted in power system stability analysis, where smaller angles indicate more stable stable solutions~\cite{kundur1994power,vournas1996relationships}. 

With these choices, the complex voltage at a bus is uniquely determined by the complex voltages of its neighbors, which leads to a natural consistency condition for a solution. Given $\bd v$, let $f$ be a function that takes $\bd v$ and performs the circle intersection operation (choosing a unique solution as described in the last paragraph). Then a vector $\bd v$ is a solution to the power flow problem if and only if $\bd v=f(\bd v)$. That is, $\bd v$ is a fixed point of $f$. Note that if two circles do not intersect at a bus, then we can declare that $\bd v$ is not a fixed point.

Here, we use the three bus network in Fig.~\ref{fig:fixedpoint} to illustrate an algorithm to solve the power flow problem. The line admittance of all the branches are $1-j\cdot1.5$. Bus $1$ is considered to be a slack bus with a voltage of $1$ p.u., while buses $2$ and $3$ are considered to be PQ buses. Initially, the voltage $\mathit{v}_{2}=v_{2,r}+jv_{2,i}$ at bus $2$ is fixed with an initial guess. Based on $\mathit{v}_{2}$, the real and reactive power circles at bus $3$ can be calculated. If these circles intersect with each other, the one with the higher voltage magnitude would be assigned as the value for $v_3$. Then, the voltage at bus $3$ is fixed and intersections of the two circles at bus $2$ are used to update $v_2$. This is repeated until the convergence is achieved. Tap changing transformers are modeled with fixed tap ratios and incorporated into the admittance matrix using $\pi$ equivalent representations. Next, we describe the algorithm for a general network.

%% file: algorithm.tex
%% !TEX root=main.tex
\subsection{Description of the Algorithm}
\begin{figure*}[ht!]
  \centering
  \begin{subfigure}{0.31\textwidth}
    \includegraphics[width=\textwidth]{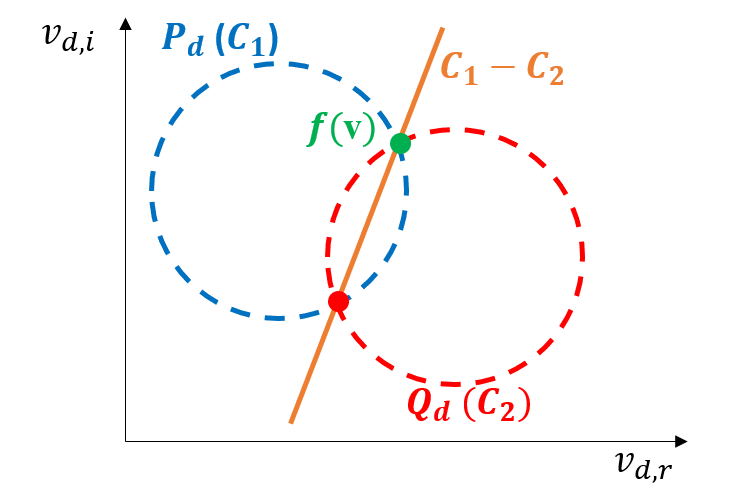}
    \caption{Equation of a line passing through common points of two circles.}
    \label{fig:rep2}
  \end{subfigure}
  ~
  \begin{subfigure}{0.31\textwidth}
    \includegraphics[width=\textwidth]{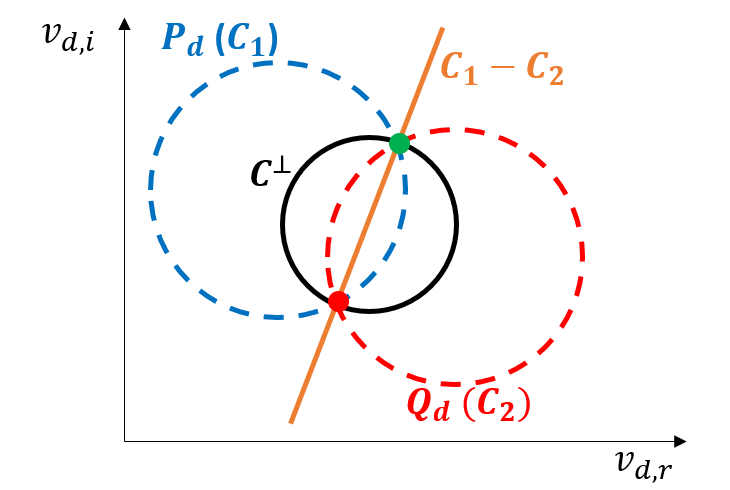}
    \caption{Orthogonal circle passing through common points of two circles.}
    \label{fig:rep3}
  \end{subfigure}
  ~~~~~~
  \begin{subfigure}{0.31\textwidth}
    \vspace{4mm}
    \includegraphics[width=\textwidth]{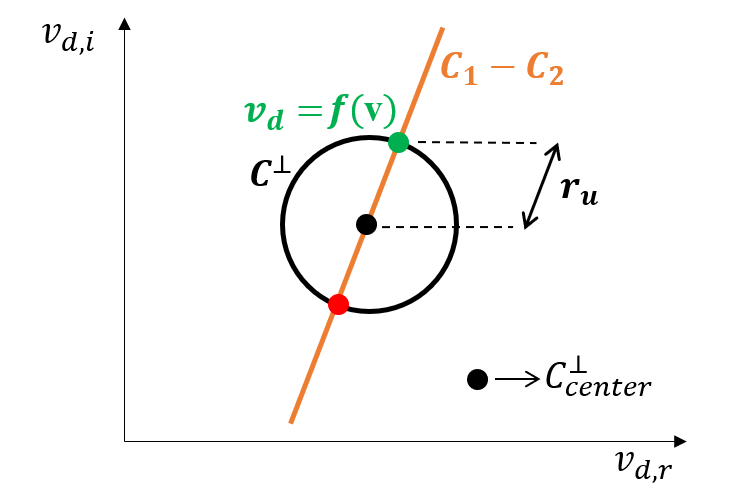}
    %		\vspace{0.3mm}
    \caption{Point of intersection of orthogonal circle and the line passing through common points is $\mathit{v}_{d} = \mathit{f(\bd v)}$.}
    \label{fig:rep4}
  \end{subfigure}
  \caption{Geometrical illustration of calculating the voltage solution at a bus.}
  \label{fig:overrall_rep}
\end{figure*}
For an n-bus system, to start the algorithm, the voltages at all the buses in the system are fixed with an initial guess. Then the voltage solution at a bus is updated using its neighbors. This is repeated for all buses, which we call a round of the algorithm. The algorithm terminates if none of the buses update their complex power in a round or when the complex power mismatch is less than the tolerance set by the user. Algorithm~\ref{alg:fixedpoint} presents the pseudo code for a system with only PQ buses. For a system with mixed PQ and PV buses, a similar algorithm can be used, which is given in the Appendix~\ref{subsec:handle_PQ_PV_bus}.

In our implementation of Algorithm~\ref{alg:fixedpoint}, we sweep through all of the buses in one round, as illustrated in Fig.~\ref{fig:update} using a lexicographical order. The exact algorithm to find the intersection will be explained in more details in Section~\ref{sec:FPE_powerflow_prob}. The exact order of updates is not constrained by the algorithm, although it is an interesting question to see if there exist an ``optimal'' update order in some sense.

\begin{algorithm}[!h]
\caption{Fixed point algorithm for system with only PQ buses.}\label{alg:fixedpoint}
\begin{algorithmic}[1]
  %  \Procedure{myalgo}{}
  \Statex \textbf{Input~~:} $\mathit{P}_i$, $\mathit{Q}_i$ for bus $i=2,\cdots,n$,
  \Statex ~~~~~~~~~~~Tolerance $\delta$ for the stopping criterion.
  \Statex \textbf{Output:} $\mathit{v}_i$ for bus $i=2,\cdots,n$.
  %\vspace{-1em}
  \State Initialize voltages at all buses, $\mathit{v}_i$ for $i=2,\cdots,n$;
  \State  Let the neighboring bus index be $k$;
  \Statex where $k$ $\in$ $\mathbb{N}_m$.
  \Comment{$\mathbb{N}_m$}{: Neighboring buses of bus $m$.}
  \State Calculate the power mismatch ($\Delta\mathbb{S}$);

  \While{($\Delta\mathbb{S}) > \delta$}
  \Comment{}{Convergence criteria.}
  \For{$m = 2$, $\cdots$, $n$}
  \LState Calculate $(\mathit{c}_{p},\mathit{r}_{p})$  at bus $m$, $\forall$ $k$ $\in$ $\mathbb{N}_m$;
  \LState Calculate $(\mathit{c}_{q},\mathit{r}_{q})$  at bus $m$, $\forall$ $k$ $\in$ $\mathbb{N}_m$;
  \LState Calculate the voltage $\mathit{v}_m$ for bus $m$;
  \LState $\mathit{v}_m = \mathit{v}_m$;
  \Comment{}{Update current bus ($m$) voltage.}
  \EndFor
  \LState Calculate ($\Delta\mathbb{S}$);
  \EndWhile
  \State \Return $\mathit{v}_m$ $\forall$ buses $m=2,\cdots,n$;
  %  \EndProcedure
\end{algorithmic}
\end{algorithm}

It is possible that the circles do not intersect at a bus either at the start of the algorithm or during one of the iterations. In these cases, we simply restart the algorithm with a new initial guess. We note that for feasible problems with PQ buses, we have never observed the non-intersection of the circles. However for non-feasible problems, it is observed that the algorithm will have non-intersection of circles and further updates in the iterative process are not possible. For mixed PQ and PV systems, it is possible for bad initial guesses to lead to non-intersection behaviors, especially when the loading is very heavy. We will provide more details in Section~\ref{sec:results}.

The proposal algorithm is similar in spirit to the Gauss-Seidel updates used in solving linear equations, where computed results are used as soon as they become available~\cite{hageman2012applied}. These type of algorithms are very efficient in memory but can converge slowly. In Section~\ref{sec:results}, we show that our algorithm does not suffer from these potential issues and are competitive with Newton-type methods for large systems.   
% In cases where the two circles \eqref{math_p} and \eqref{math_q} at bus $d$ has no intersection then we restart with random initialization value. The scenario for \eqref{math_p} and \eqref{math_q} to not have an intersection is almost never seen in a PQ bus system. However, this non-intersection behavior is observed in a very few cases of PQ and PV bus system especially when the power system is heavily loaded as compared to its base case loading.

%% file: fixed_point_equation.tex
%% !TEX root=main.tex

In Algorithm~\ref{alg:fixedpoint}, the main computation step is to find the intersection of two circles. At a first glance, this operation is almost trivial and there are many different ways to compute the intersections (see, e.g.~\cite{weisstein2003circle}). However, the numerical implementation of an intersection algorithm can experience subtle but critical issues. First of all, this operation is called upon many times in the algorithm, and small errors can propagate and result in slower convergence speeds. Second of all, the circles can have very small or large radii.  For example, if lines are close to being purely inductive (high $X/R$ ratio), then the reactive circle becomes a circle with a very large radius and straightforward algorithms would run into numerical instabilities. Finally, since finding the intersections takes most of the time in Algorithm~\ref{alg:fixedpoint}, it would be desirable to get as close to a closed-form solution as possible. Therefore, we use an unconventional representation of circles developed by~\cite{middle88} to provide a robust and efficient algorithm to find the intersection of circles. For ease of exposition, we focus on system with PQ bus. Analogous results can be derived for PV buses.

Fig.~\ref{fig:overrall_rep} outlines the steps we take to find the intersection of the active and reactive power circles. First, we find the line through the two circles (Fig.~\ref{fig:rep2}). Then, we find the smallest circle (called the orthogonal circle) that passes through the intersecting points of the original circles (Fig.~\ref{fig:rep3}). Next, we find the intersection of the line with the orthogonal circle (Fig.~\ref{fig:rep4}). It turns out that if we view the circles as vectors in a vector space, the above computations can be thought as vector manipulations, which is simple to perform and numerically stable. In the rest of this section, we develop this theory based on the material in~\cite{middle88}.

\begin{remark}
    Note that we do not necessarily require the step in Fig.~\ref{fig:rep3}. The key reason to compute another circle is that the numerical accuracy of the algorithm improves if the smallest possible circle is used to compute the intersection points as pointed in~\cite{middle88}. Therefore, we use the orthogonal circle, which is the smallest possible circle that contains the intersection points in our algorithm.
    % improves The fixed-point representation $\mathit{f(\bd v)}$ can also be derived as an intersection between the common line ($\mathit{C}_{1}$ - $\mathit{C}_{2}$) and either $\mathit{C}_{1}$ or $\mathit{C}_{2}$. In fact, the use of a smaller circle in between $\mathit{C}_{1}$ and $\mathit{C}_{2}$ gives a more accurate result. However, \cite{middle88} shows that the intersection between the common line and orthogonal circle gives even a better result since the common line crosses the $\mathit{C}^{\perp}$ at a less acute angle. Hence, in this paper, the fixed-point representation is obtained by the intersection of the common line and orthogonal circle.
\end{remark}

\subsection{3-Tuple Vector Representation of Circles}
\label{sec:subsec_3tup_rep}
Instead of the traditional center/radius parameterization, we can describe all of the points $\bd x \in \mathbb{R}^2$ on a circle by the following equation:
\begin{align}
   \mathit{a} (\mathit{{\textbf{x}}}\cdot \mathit{{\textbf{x}}}) + \mathit{{\textbf{b}}}\cdot\mathit{{\textbf{x}}}+\mathit{c} = 0, \label{eq:base_form}
\end{align}
where $\cdot$ denotes the dot product between two vectors. The form in \eqref{eq:base_form} allows us to describe a circle using a three tuple $(a,\bd b,c)$. Note that this presentation is not unique, since scaling all of the parameters by a scalar does not change the points that satisfy \eqref{eq:base_form}. If $a$ is not zero, we will scale parameters such that $a=1$. In this notation, the circles described by the real and reactive power equations \eqref{math_p} and \eqref{math_q} becomes
\begin{equation}
\left( \mathit{a}_{p}, \mathit{{\textbf{b}}}_{p}, \mathit{c}_{p} \right) = \left(1, \begin{bmatrix}	\frac{t_{d,2}}{t_{d,1}} & \frac{t_{d,3}}{t_{d,1}} \end{bmatrix}^{T}, -\frac{\mathit{p}_{d}}{t_{d,1}} \right),\label{eq:p_tup}
\end{equation}
and
\begin{equation}
\left( \mathit{a}_{q}, \mathit{{\textbf{b}}}_{q}, \mathit{c}_{q} \right) = \left( 1, \begin{bmatrix} -\frac{t_{d,3}}{t_{d,4}} &  \frac{t_{d,2}}{t_{d,4}}\end{bmatrix}^{T}, -\frac{\mathit{q}_{d}}{t_{d,4}} \right).
\label{eq:q_tup}
\end{equation}
In these representations, the circles shift gracefully and the same calculations can be applied to a wide range of parameter values even if they approach zero or infinity.
% For example, if the lines are inductive, $t_{d,1}=0$ for all $d$, and \eqref{eq:p_tup} can still be used without numerical issues (unlike the center-radius approach in~\cite{weisstein2003circle}.)

Next, we separate the fixed parameters in the system (e.g., admittance values) and the voltages. Given a bus $d$, let $d_1,d_2,\cdots,d_k$ its neighboring nodes. Let $\bd g_d=\begin{bmatrix} g_{d_1,d} & g_{d_2,d} & \cdots & g_{d_k,d} \end{bmatrix}$ denote the the vector of conductances between bus $d$ and its neighbors. Similarly, let $\bd b_d=\begin{bmatrix} b_{d_1,d} & b_{d_2,d} & \dots & b_{d_k,d} \end{bmatrix}$ denote the vector of susceptances. Let $\bm{1}$ denote the vector of all $1$'s of the appropriate length. To represent the voltages of the neighboring buses, we use a vector $\bd u$ formed by concatenating the real and imaginary voltages:
\begin{equation*}
\bd u=\begin{bmatrix}\mathit{v}_{d_1,r} & \mathit{v}_{d_2,r} & \cdots & \mathit{v}_{d_k,r} & \mathit{v}_{d_1,i} & \mathit{v}_{d_2,i} & \cdots & \mathit{v}_{d_k,i} \end{bmatrix}^T.
\end{equation*}
Then, we can rewrite \eqref{eq:p_tup} and \eqref{eq:q_tup} as
\begin{align}
\left( \mathit{a}_{p}, \mathit{{\textbf{b}}}_{p}, \mathit{c}_{p} \right) = \left( 1,
\begin{bmatrix}
	-\bm{\alpha} & \bm{\delta} \\ -\bm{\delta} & -\bm{\alpha}
\end{bmatrix} \bd{u},\dfrac{\mathit{p}_{d}}{\bm{1}\cdot \bd{g}_d} \right), \label{eq:3tuple_P}
\end{align}

\begin{align}
\left( \mathit{a}_{q}, \mathit{{\textbf{b}}}_{q}, \mathit{c}_{q} \right) = \left( 1,
\begin{bmatrix}
	-\bm{\beta} & -\bm{\gamma} \\ \bm{\gamma} & -\bm{\beta}
\end{bmatrix} \bd u,\dfrac{\mathit{q}_{d}}{\bm{1} \cdot \bd{b}_d} \right), \label{eq:3tuple_Q}
\end{align}
where
\begin{align}
\alpha = \dfrac{\bd g_{d}}{\bm{1}\cdot \bd{g}_d}, \beta = \dfrac{\bd b_d}{\bm{1} \cdot \bd{b}_d}, \gamma = \dfrac{\bd g_{d}}{\bm{1} \cdot \bd{b}_d}, \delta = \dfrac{\bd b_d}{\bm{1}\cdot \bd{g}_d}. \notag
\end{align}
% \begin{align}
% \mathit{\textbf{g}}_{k,d} &= \begin{bmatrix}\mathit{g}_{1,d} & \mathit{g}_{2,d} & \cdots & \mathit{g}_{d-1,d} & \mathit{g}_{d+1,d}\end{bmatrix}, \notag \\
% \mathit{\textbf{b}}_{k,d} &= \begin{bmatrix}\mathit{b}_{1,d} & \mathit{b}_{2,d} & \cdots & \mathit{b}_{d-1,d}& \mathit{b}_{d+1,d}\end{bmatrix}, \notag \\
% \bd{u}_{1} &= \begin{bmatrix}\mathit{v}_{1,r} & \mathit{v}_{2,r} & \cdots & \mathit{v}_{d-1,r}& \mathit{v}_{d+1,r}\end{bmatrix}^T, \notag \\
% \bd{u}_{2} &= \begin{bmatrix}\mathit{v}_{1,i} & \mathit{v}_{2,i} & \cdots & \mathit{v}_{d-1,i} & \mathit{v}_{d+1,i}~ \end{bmatrix}^T, \notag
% \end{align}

%  the Using the three tuple notation in \eqref{eq:base_form} also allows us to apply vector calculations. \eqref{eq:p_tup} and \eqref{eq:q_tup} can be written a vectorized form as shown in \eqref{eq:3tuple_P} and \eqref{eq:3tuple_Q} respectively. The simplification steps to obtain the vector form is presented in Appendix~\ref{subsec:vector_simplifi}. This vectorized representation enables ease of calculation in terms of speed and complexity in derivations.
%
%
%
% $\bd{u}$ represents the neighboring voltages of bus ``d" such that $\bd{u} = \begin{bmatrix} \bd{u}_{1} & \bd{u}_{2} \end{bmatrix}^T$ and buses 1, 2, $\cdots$, $d-1$ and $d+1$ are the neighboring buses to bus $d$.

If needed, the centers and radii of power flow circles can be computed easily from \eqref{eq:3tuple_P} and \eqref{eq:3tuple_Q}:
\begin{align}
&\mathbb{\textbf{Center}}_{p} = \dfrac{-\mathit{\textbf{b}}_{p}}{2}, \mathbb{\textbf{Center}}_{q} = \dfrac{-\mathit{\textbf{b}}_{q}}{2}, \label{eq:cpq_tuple} \\
&\mathit{r}^2_{p} = \left(\dfrac{\mathit{\textbf{b}}_{p}\cdot \mathit{\textbf{b}}_{p}}{4} - \mathit{c}_{p}\right), \mathit{r}^2_{q} = \left(\dfrac{\mathit{\textbf{b}}_{q}\cdot \mathit{\textbf{b}}_{q}}{4} - \mathit{c}_{q}\right),\label{eq:rpq_tuple}
\end{align}
where $\mathbb{\textbf{Center}}_{p}$ and $\mathbb{\textbf{Center}}_{q}$ are the centers of the real and reactive power circles, $r_p$ and $r_q$ are the radii, respectively.

\subsection{Line Passing Through Intersection Points of the Power Flow Circles}
	Given two circles $C_{1} = (1,\mathit{{\textbf{b}}}_{1},\mathit{c}_{1} )$ and $\mathit{C}_{2} = ( 1,\mathit{{\textbf{b}}}_{2},\mathit{c}_{2} )$, the line passing through their points of intersection is described by $C_1-C_2$, provided the circles intersect. More formally, $C_1-C_2$ is
	\begin{align}
	\mathit{C}_{1}-\mathit{C}_{2}= ( 0,\mathit{{\textbf{b}}}_{1}-\mathit{{\textbf{b}}}_{2},\mathit{c}_{1}-\mathit{c}_{2} )=(0,\bd L_2,L_3) \label{eq:line_eq_gen}
	\end{align}
and describes the points $\bd v_d$ that satisfies the equation
	\begin{equation}
	\mathit{{\textbf{L}}}_{2} \cdot \mathit{{\textbf{v}}}_{d} + \mathit{L}_{3}=0, \label{eq:line_form}
	\end{equation}
	where
	\begin{equation} \mathit{{\textbf{v}}}_{d} = \begin{bmatrix} \mathit{v}_{d,r} \\ \mathit{v}_{d,i} \end{bmatrix}.\notag \end{equation}
Substituting \eqref{eq:3tuple_P} and \eqref{eq:3tuple_Q} into \eqref{eq:line_form}, we have the line described by
\begin{equation}
	\left ( 0,\begin{bmatrix}-\bm{\alpha}+\bm{\beta} & \bm{\gamma}+\bm{\delta} \\ -\left(\bm{\gamma}+\bm{\delta}\right) & -\bm{\alpha}+\bm{\beta} \end{bmatrix} \bd u,\dfrac{\mathit{p}_{d}}{\bm{1} \cdot \bd{g}_d} + \dfrac{\mathit{q}_{d}}{\bm{1} \cdot \bd{b}_d} \right). \label{eq:common_line}
	\end{equation}
	\subsection{Orthogonal Circle}
  In principle, we can use the line computed in \eqref{eq:common_line} to find the intersection points by intersecting that line with one of the active or reactive circles. However, the numerical accuracy and stability can suffer because the line may intersect the circles at a very acute angle. Therefore, it is more desirable to use the orthogonal circle for calculations. Geometrically, the orthogonal circle is the smallest circle that passes through the two intersection points. Algebraically, we label it as $C^{\perp}$. Again, the parameters of this circle can be computed from \eqref{eq:3tuple_P} and \eqref{eq:3tuple_Q} via simple algebra~\cite{middle88,baker94}:\footnote{The original formula given in~\cite{middle88} is in fact incorrect and the right formula is given in~\cite{baker94}.}
	\begin{align}
	\mathit{C}^{\perp}	&= \left ( a^{\perp}, \bd b^{\perp}, c^{\perp}\right) \notag \\
  &=	\left ( 1,\dfrac{\mathit{{\textbf{b}}}_{1} + \mathit{{\textbf{b}}}_{2}}{2} + \dfrac{\left( \mathit{{\textbf{b}}}_{2} - \mathit{{\textbf{b}}}_{1} \right)  \left( \mathit{k}^{2}_{1} - \mathit{k}^{2}_{2}\right)}{2  \left\| \mathit{{\textbf{b}}}_{1} - \mathit{{\textbf{b}}}_{2}\right\|^{2}} ,\right. \notag \\
	&\left.~~~~\dfrac{\mathit{c}_{1} + \mathit{c}_{2}}{2} + \dfrac{\left( \mathit{c}_{2} - \mathit{c}_{1} \right)  \left( \mathit{k}^{2}_{1} - \mathit{k}^{2}_{2}\right)}{2  \left\| \mathit{{\textbf{b}}}_{1} - \mathit{{\textbf{b}}}_{2}\right\|^{2}} \right),\label{eq:c_ortho_gen}
	\end{align}
	where
	\begin{align*}
	\mathit{k}^{2}_{1} &= \|\mathit{{\textbf{b}}}_{1}\|^2 - 4\mathit{a}_{1}\mathit{c}_{1}, \\
	\mathit{k}^{2}_{2} &= \|\mathit{{\textbf{b}}}_{2}\|^2 - 4  \mathit{a}_{2}  \mathit{c}_{2}.
	\end{align*}
	Here, $\|~\|$ is the standard $l_2$ norm. The center and the radius of the orthogonal circle is given by
	\begin{equation}
	\mathit{\textbf{Center}}^{\perp} = \dfrac{-\bd b^{\perp}}{2}
	\end{equation}
  and
	\begin{equation}
	\mathit{r}_{u} = \sqrt{\dfrac{\bd b^{\perp} \cdot \bd b^{\perp}}{4} - c^{\perp}}. \label{eq:c_ortho_radius}
	\end{equation}

Substituting \eqref{eq:3tuple_P} and \eqref{eq:3tuple_Q} in \eqref{eq:c_ortho_gen}, we get
	\begin{align}
	\mathit{C}^{\perp} &= \left ( 1, \frac{1}{2} \bd{M}_{B} \bd u + \dfrac{\left( \dfrac{\left\|\bd{u}\right\|^2}{2} \cdot \mathit{K}_{c} - 2 \cdot l\right)  \bd{M}_{A} \bd u }{\left \|\bd{M}_{A} \bd{u} \right \|^2},\right.\notag \\
	&\left. \dfrac{1}{2} \left( \dfrac{\mathit{p}_{d}}{\bm{1} \cdot \bd{g}_d} - \dfrac{\mathit{q}_{d}}{\bm{1}\cdot\bd{b}_d}\right) +\dfrac{l \left( 2 l - \dfrac{\left\|\bd{u}\right\|^2}{2} \mathit{K}_{c}\right)}{\left \| \bd{M}_{A} \bd{u} \right \|^2}\right), \label{eq:c_orth}
	\end{align}
	where
	\begin{align}
	\bd{M}_{A} &= \begin{bmatrix} \bm{\alpha}-\bm{\beta} & -\left(\bm{\gamma}+\bm{\delta}\right) \\ \bm{\gamma}+\bm{\delta} & \bm{\alpha}-\bm{\beta} \end{bmatrix}
	\bd{M}_{B}= \begin{bmatrix} -\left(\bm{\alpha}+\bm{\beta}\right) & \bm{\delta}-\bm{\gamma} \\ -\left(\bm{\delta}-\bm{\gamma}\right) & -\left(\bm{\alpha}+\bm{\beta}\right) \end{bmatrix}, \notag \\
  l&=\dfrac{\mathit{p}_{d}}{\bm{1} \cdot \bd{g}_d} + \dfrac{\mathit{q}_{d}}{\bm{1} \cdot \bd{b}_d},~\mathit{K}_{c} = \left( \|\bm{\alpha}\|^2 + \|\bm{\delta}\|^2 \right) - \left( \|\bm{\gamma}\|^2 + \|\bm{\beta}\|^2\right). \notag
	\end{align}

	\subsection{Point of Intersection}
	Next, we find the point of intersection (Fig.~\ref{fig:rep4}). These points are found at a distance of $\mathit{r}_{u}$ from the center of orthogonal circle along the line computed in~\eqref{eq:common_line}. Through simple algebra, we compute the point of intersection, that is, the updated voltage at bus $d$ given by
   % Let there be a point $\mathit{\textbf{X}}_{t}$ which is $\mathit{t}$ distance away from another point $\mathit{\textbf{X}}_{0}$ on the same straight line of the form $\mathit{\textbf{a}} \cdot \mathit{\textbf{X}} + \mathit{b} = 0$. Then the parametric representation of the point $\mathit{\textbf{X}}_{t}$ is given by \eqref{eq:param_line_gen}. Similarly, the parametric form for the common point $\mathit{f(\bd v)}$ in Fig.~\ref{fig:rep4} is given by \eqref{eq:param_line_gen1}.
	\begin{align}
	\begin{bmatrix} v_{d,r} \\ v_{d,i} \end{bmatrix} &= \mathit{\textbf{Center}}^{\perp} \pm \mathit{r}_{u}  \dfrac{\bd{R} \mathit{{\textbf{L}}}_{2}}{\left \| \mathit{{\textbf{L}}}_{2}\right \|^{2}} \nonumber \\
  &= -\frac{\bd b^{\perp}}{2} \pm \sqrt{\dfrac{\bd b^{\perp} \cdot \bd b^{\perp}}{4} - c^{\perp}} \cdot \dfrac{\bd{R} \mathit{{\textbf{L}}}_{2}}{\left \| \mathit{{\textbf{L}}}_{2}\right \|^{2}}
  , \label{eq:param_line_gen1} \\
  &= \frac{-1}{4} \bd{M}_{B} \bd u - \dfrac{\left( \dfrac{\left\|\bd{u}\right\|^2}{2} \cdot \mathit{K}_{c} - 2 \cdot l\right)  \bd{M}_{A} \bd u }{2 \cdot \left \|\bd{M}_{A} \bd{u} \right \|^2} \pm \notag \\ 
  & \left(\dfrac{1}{4} \cdot \left\| \frac{1}{2} \bd{M}_{B} \bd u + \dfrac{\left( \dfrac{\left\|\bd{u}\right\|^2}{2} \cdot \mathit{K}_{c} - 2 \cdot l\right)  \bd{M}_{A} \bd u }{\left \|\bd{M}_{A} \bd{u} \right \|^2} \right\|^2 - \right. \notag \\
  & \left.\dfrac{1}{2} \left( \dfrac{\mathit{p}_{d}}{\bm{1} \cdot \bd{g}_d} - \dfrac{\mathit{q}_{d}}{\bm{1}\cdot\bd{b}_d}\right) +\dfrac{l \left( 2 l - \dfrac{\left\|\bd{u}\right\|^2}{2} \mathit{K}_{c}\right)}{\left \| \bd{M}_{A} \bd{u} \right \|^2}\right)^{\dfrac{1}{2}} \cdot \dfrac{\bd{R} \mathit{{\textbf{L}}}_{2}}{\left \| \mathit{{\textbf{L}}}_{2}\right \|^{2}}. \notag
	\end{align}
	where
	\begin{equation}
	\bd R = \begin{bmatrix} 0 & -1 \\ 1 & 0 \end{bmatrix} \notag
	\end{equation}
and $\bd b^{\perp}$ is given by \eqref{eq:c_ortho_gen} and $\bd L_2$ is given by \eqref{eq:line_eq_gen}. To choose one solution or a sign in \eqref{eq:param_line_gen1}, we will pick the one that leads to the higher voltage magnitude. Note that in \eqref{eq:param_line_gen1}, most of the computation can be done offline since they only involve the admittance parameters. Applying \eqref{eq:param_line_gen1} to every bus $d$ also gives us an analytical form of the fixed point equation for the complex voltages. 

%% file: numerical.tex
%% !TEX root=main.tex
%		The criteria for convergence of the algorithm is, the sum of power mismatches at all the buses is equal to \num{1e-3} p.u and base 100 MVA.

    In this section, simulation studies on standard IEEE test cases, specifically with the 4, 14, 30, 39, 57, 118, 2383 and 3375 bus systems. Their information are obtained from the Matpower software~\cite{ZimmermanEtAl2011}. The proposed fixed point algorithm are compared with other power flow methods at nominal loading and heavy loading conditions. We also test the sensitivity of these algorithm to the initializaton points. In particular, we will consider five algorithms: 1) \textbf{FP}, our proposed algorithm; 2) \textbf{GS}, the standard Gauss-Seidel algorithm; 3) \textbf{NR}, the standard Newton-Raphson algorithm; 4) \textbf{FDLF}, the fast decoupled load flow algorithm and 5) \textbf{Iwamoto}, a Jacobian-based adjustable step size method (sometimes called non-divergent power flow method)~\cite{Iwamoto81}.
    % The comparative studies include heavy loading, and effect of random initialization.

\subsection{Performance of Proposed Method}
\begin{figure}[!h]
		\centerline{\includegraphics[width=\linewidth]{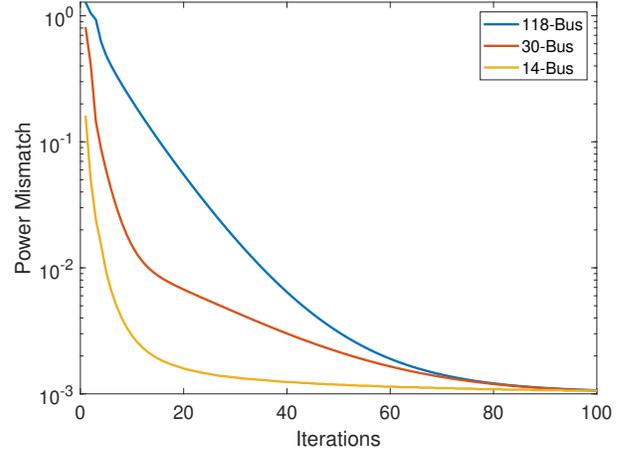}}
		\caption{Semi-log plot of the convergence of the fixed point algorithm for IEEE standard systems at base case loading.}
		\label{fig:con_RFPPF_IEEE_cases}
	\end{figure}

First we study the convergence of the proposed FP algorithm on the standard $14$-, $30$- and $118$-IEEE bus systems are shown in Fig.~\ref{fig:con_RFPPF_IEEE_cases} under nominal loading conditions. More precisely, we take the data from MatPower~\cite{ZimmermanEtAl2011} and run our fixed point algorithm to test its convergence. As shown in Fig.~\ref{fig:con_RFPPF_IEEE_cases}, for these standard cases, the fixed point algorithm converges in tens of iterations. Since each iteration is cheap to compute, the convergence time is in 10s of milliseconds. %\textcolor{violet}{We note that the time it takes to finish one iteration is much faster than say a conventional Newton type algorithm, where the inverse of the Jacobian need to be computed.} \textcolor{blue}{This previous statement made the reviewers think that we are trying to say that RPF method is faster than NR method and they asked us to compare the speed. I think we have to admit that our method is obviously slower than NR method. So, maybe we can remove the violet colored text above or is it okay to keep it?}

\begin{figure}[H]
		\centerline{\includegraphics[width=\linewidth]{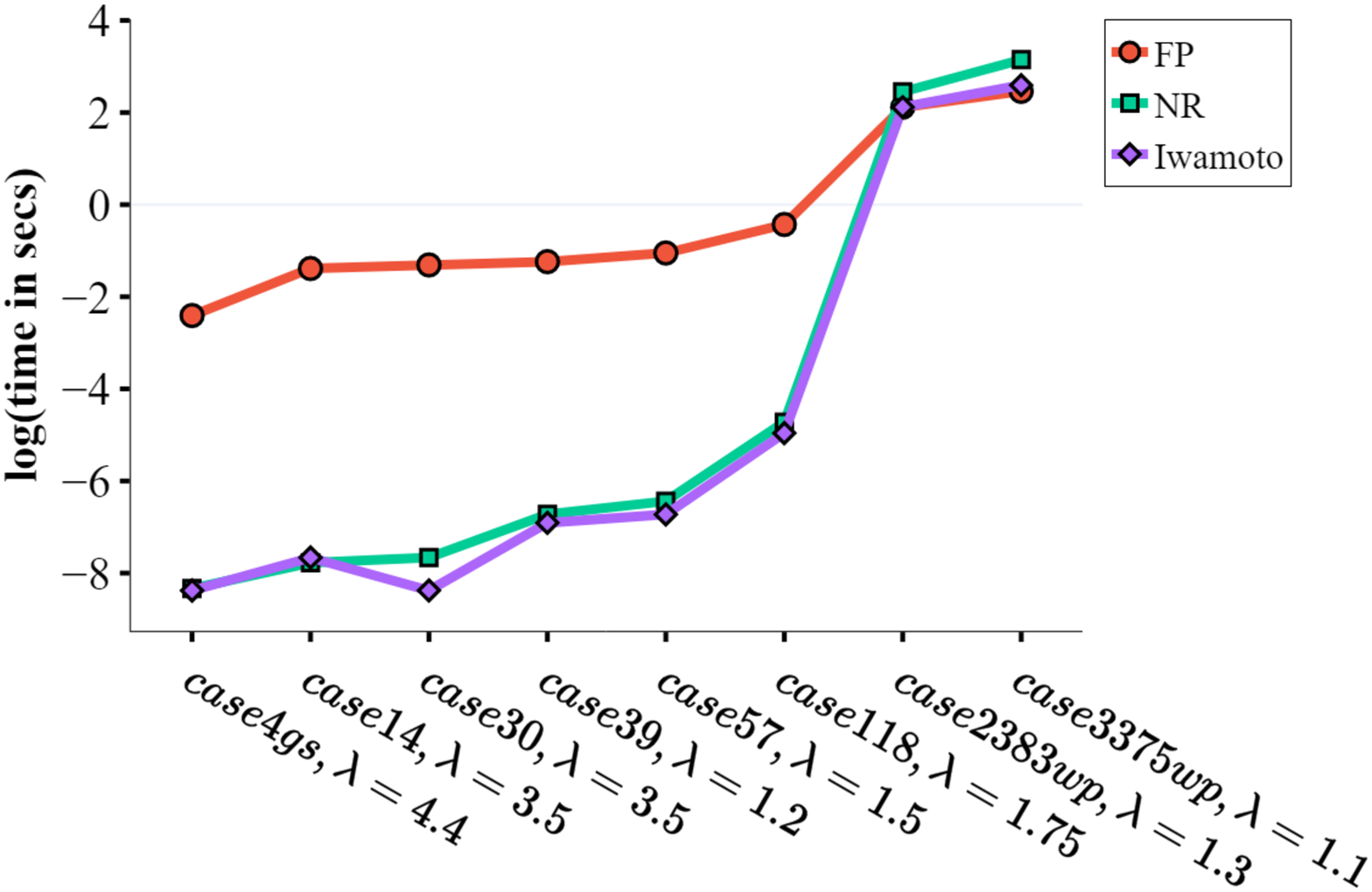}}
		\caption{Time taken to attain the desired precision by FP, NR and Iwamoto's method for different size systems under different load scaling parameters $\lambda$. As the systems become large, FP becomes competitive with the other two methods.}
		\label{fig:time_comp}
	\end{figure}

Next, we compare the convergence speed of FP with NR and Iwamoto algorithms for a variety of systems. In addition to using the nominal loading, we introduce a scaling parameter $\lambda$ to scale the loads and generations by multiplying both the active and reactive powers by $\lambda$. The performance of all three algorithms are shown in Fig.~\ref{fig:time_comp} in a log-scale. As expected, when the networks are small, NR and Iwamoto methods converges very quickly since they only require a few (sometimes only 1) updates. As the test cases becomes bigger, the FP method start to catch up. For the largest two test systems (2383-bus and 3375-bus), FP is comparable with Iwamoto and faster than NR. In addition, convergence times for FP is much less sensitive to the system size, which is expected since there are not matrix computations. 

    \begin{figure}[H]
		\centerline{\includegraphics[width=\linewidth]{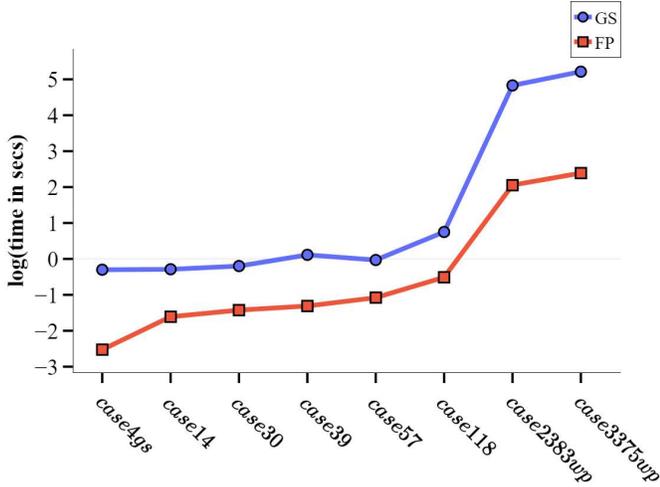}}
		\caption{Time taken to attain the desired precision by GS and FP method for different size systems at base case loading. FP is faster than GS for all of the IEEE systems.}
		\label{fig:time_comp2}
	\end{figure}
	
Fig.~\ref{fig:time_comp2} compares the convergence speed of FP and GS for a variety of systems. GS calculates the voltage solutions at every bus in the system via a lexicographical approach. Even though GS is partly similar to the proposed approach, the fixed point equation used to calculate the voltages in both the methods are different. This difference makes the FP perform faster when compared to GS as shown in Fig.~\ref{fig:time_comp2}.

\subsection{Heavily Loaded Networks}
    \begin{figure}[ht]
		\centerline{\includegraphics[width=\linewidth]{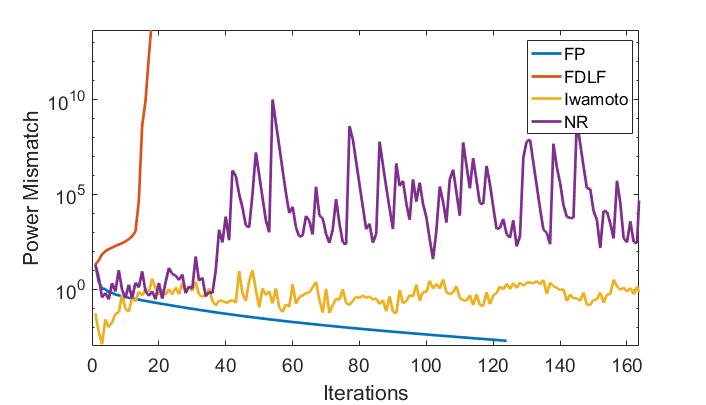}}
		\caption{Convergence performance of FP compared with NR, Iwamoto and FDLF. The test system is the IEEE 14-bus system with loads scaled by a factor of 3.99. At this loading, the FP method converges reliably to the solution while the other methods either diverges or exhibit oscillatory behavior.}
		\label{fig:RFPPF_heavy_other_methods}
	\end{figure}
	
The more challenging setting, and the setting the FP algorithm is designed to address, is when the systems are heavily loaded. For example, we take the $14$-bus network and increase all of the loads by a factor of 3.99. This loading is still feasible, but is very close to the loadability limit of the system.  Figure~\ref{fig:RFPPF_heavy_other_methods} presents the convergence comparison of the methods. As expected, NR diverges~\cite{VanCutsemEtAl2007,VanCutsem2000} since the Jocabian matrix is very ill-conditioned around the solution. The FDLF also diverges, while it does not face conditioning problems since the Jacobian is approximated by a fixed matrix in decoupled load flow, the update direction provided by the fixed Jacobian becomes invalid and the algorithm diverges very quickly. Interestingly, the optimal multiplier method (Iwamoto~\cite{Iwamoto81})  also becomes unstable because of numerical issues. More precisely, the multiplier $\mu$ is used to control the step size in the following update:
\begin{equation}
 \Delta\mathit{V} = \mu\cdot\left(J^{-1}\right)\cdot\left(\Delta\mathit{S}\right). \label{eq:opt_mul}
\end{equation}
However, if $J$ is very ill-conditioned, then we are essentially multiplying a very small number by a very large number, which creates problems due to finite machine precision.

Similar behaviors to Fig.~\ref{fig:RFPPF_heavy_other_methods} are observed in  IEEE $4$-, $30$- and $118$-bus systems for load multiplier of 4.5, 3.65 and 1.78, respectively. In contrast, our proposed method is able to converge even under these conditions since it does not use the power flow Jacobian. This illustrates the envisioned utility of the proposed FP method in practice. An operator can use conventional power flow solvers and when they do not converge, instead of fine tuning parameters or trying many different initialization points, the FP algorithm can be used as a viable tool to obtain convergence.

% issues when inverting the matrix. This issue is also studied in \cite{Iwamoto81} in terms of single and double precision. \textcolor{blue}{Similar observations are made for IEEE $4$-, $30$- and $118$-bus systems for load multiplier of 4.5, 3.65 and 1.78 respectively.} In our simulation, we use \eqref{eq:opt_mul} to compute the solution of an optimal multiplier based method. Even though the optimal multiplier $\mu$ changes the step size of the Jacobian, the Jacobian becomes highly sensitive when nearing to singularity and an incorrect $\Delta\mathit{V}$ is predicted. This is verified by calculating the condition number of the Jacobian matrix under the loading conditions. 

% In contrast, our proposed method is able to converge even under these conditions since it does not use the Jacobian. Similar behavior is observed for other networks.

% It is found to be $1.68\cdot10^{3}$, a large number indicating that the precision is not sufficient to calculate the Jacobian that is closer to the singularity.

\begin{figure}[ht]
		\centering
	\begin{subfigure}{\columnwidth}
		\centerline{\includegraphics[width=\linewidth]{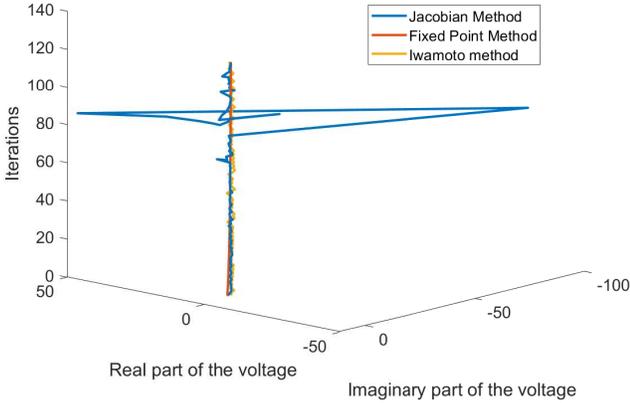}}%save_1.eps,3d_plot_v2.eps
		\caption{Updates of the real and imaginary parts of the voltage vector}
		\label{fig:subpic_1}
	\end{subfigure}
	\begin{subfigure}{\columnwidth}
		\centerline{\includegraphics[width=\linewidth]{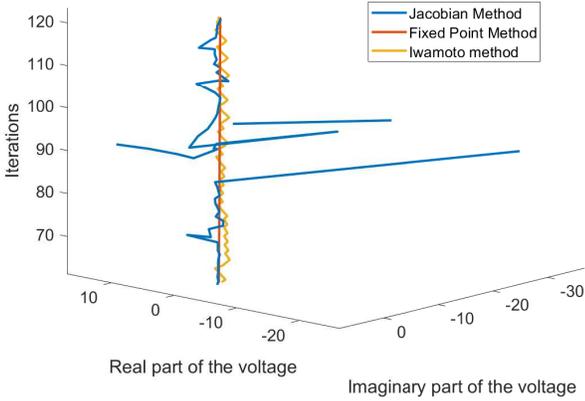}}
		\caption{Zoomed in view of Fig.~\ref{fig:subpic_1} showing the oscillatory behavior of the Iwamoto method.}
		\label{fig:subpic_2}
	\end{subfigure}
		\caption{Comparison of NR, FP and Iwamoto's voltage updates at bus $4$ in IEEE-$14$ bus system under the loading condition of Fig.~\ref{fig:RFPPF_heavy_other_methods}.}
		\label{fig:RFPPF_heavy_other_methods_comp}
	\end{figure}

Figure~\ref{fig:RFPPF_heavy_other_methods_comp} compares the performance of the FP, NR, FDLF and Iwamoto algorithms in more detail. As we can see, the NR algorithm jumps erratically in the voltage space. The Iwamoto method controls this behavior by scaling the updates by $\mu$, but even though it prevents the algorithm from diverging, it cannot converge reliably and instead oscillate around the solution. The FP algorithm again converges reliably and do not oscillate. It does not encounter the numerical instabilities of other methods since the only calculation required is the intersection of two circles (regardless of the system size), and these intersections can be handled gracefully even when the circles becomes degenerate using the algorithm in Section~\ref{sec:FPE_powerflow_prob}.

% 	\begin{figure}[!h]
% 		\centerline{\includegraphics[width=\linewidth]{new_comp_fig.eps}}
% 		\caption{\textcolor{blue}{Comparison of voltage updates at bus $4$ in IEEE-$14$ bus system under heavily loaded conditions.}}
% 		\label{fig:RFPPF_heavy_other_methods_comp}
% 	\end{figure}
	
% 	\begin{table}[H]
% 	\centering
% 	\renewcommand{\arraystretch}{1.5}
% 	\begin{tabular}{|>{\columncolor[RGB]{51, 153, 255}}c|>{\columncolor[RGB]{51, 153, 255}}c|>{\columncolor[RGB]{51, 153, 255}}c|>{\columncolor[RGB]{51, 153, 255}}c|>{\columncolor[RGB]{51, 153, 255}}c|>{\columncolor[RGB]{51, 153, 255}}c|}
% 	%https://www.w3schools.com/colors/colors_picker.asp
% 	%https://tex.stackexchange.com/questions/94799/how-do-i-color-table-columns
% 		\hline \hline
% 		Bus System  & Multiplier&NR&Adjusted step-size &GS&FP \\ \hline
% 		4-bus& 4.5| 4.54&0& 0& 0& 1 \\ \hline
% 		14-bus& 3.99| 4&0& 0& 0& 1 \\ \hline
% 		30-bus& 3.65| 3.658&0& 1& 0& 1 \\ \hline
% 		39-bus& ?|?&?& ?& ?& ? \\ \hline
% 		57-bus& ?|1.7855&?& ?& ?& ? \\ \hline
% 		2383-bus& ?|?&?& ?& ?& ? \\ \hline
% 		3210-bus& ?|?&?& ?& ?& ? \\ \hline \hline
% 	\end{tabular}
% 	\caption{\textcolor{blue}{Convergence test for various power flow algorithms with $1$ = converged and $0$ = not converged with tolerance = $0.001$.}}
% 	\label{tab:table3}
% \end{table}

\subsection{Sensitivity to Initial Conditions}
		In addition to convergence, it is important for an algorithm to be robust to the initial conditions, especially as the randomness in the system increases due to renewable integration~\cite{Costa08}. To test the performance of various algorithms to initial conditions, we take the IEEE 30-bus system at its standard loading and randomly select the starting voltages. In our experiments, we set the initial guess to be random samples from the uniform distribution on the interval $\left[1-\alpha, 1+\alpha\right]$ for various values of $\alpha$ (we always set the imaginary part to be 0) %\textcolor{red}{I ran some simulations and observed that we don't need to set the imaginary part of voltage to be zero and the results in Table.~\ref{tab:table1} still holds correct}
		, independently for each bus. Table~\ref{tab:table1} reports the number of successful convergences (defined as power mismatch convergence less than $0.001$ p.u.) for the FDLF, NR, optimal multiplier and our proposed FP methods for 100 trials.
	 % It is observed that the Algorithm~\ref{alg:fixedpoint} is insensitive to the initialization condition. This is evaluated on IEEE 30-bus system at base case loading with only PQ buses in the network. The initialization is provided by generating random samples from the uniform distribution on the interval $\left[1-\alpha, 1+\alpha\right]$ for various values of $\alpha$. For each initialization, FDXB, NR, Optimal multiplier and RFPPF methods are used and the sensitivity is calculated based on the number of successful power mismatch convergence of $0.001$ p.u. over the entire sample of 100 experiments.

%	 based on whether the converged solution is the known high voltage solution of the power system.

%	 \textcolor{blue}{Ask for tips: Should I write ``known high voltage solution" or power mismatch of $0.001$ p.u.? The former has less power mismatch than latter, which is better in paper? In reality the power mismatch I used is $0.002$ p.u., anything below this value won't describe our RFPPF as converged.}

\begin{table}[ht]
	\centering
	\renewcommand{\arraystretch}{1.5}
	\begin{tabular}{|l|c|c|c|l|}
		\hline \hline
		Initialization spread $\alpha$  & NR&FDLF&optimal multiplier&FP \\ \hline
		0.05& 100& 98& 100& 100 \\ \hline
		0.1& 64& 62& 100& 100 \\ \hline
		0.2& 4& 0& 0& 100 \\ \hline
		0.3& 0& 0& 0& 100 \\ \hline
		0.4& 0& 0& 0& 100 \\ \hline
		0.6& 0& 0& 0& 100 \\ \hline
		0.9& 0& 0& 0& 100 \\ \hline \hline
	\end{tabular}
	\caption{Convergence test of the power flow methods with random initialization for 100 trials. The initial voltages are generated identical and independently from uniform distribution of $[1-\alpha,1+\alpha]$. Our proposed method converged for every trial, where as other methods quickly stopped working once the range become moderately large.}
	\label{tab:table1}
\end{table}

As we see in Table~\ref{tab:table1}, our proposed FP method is much more robust to the value of the initial guesses than the other methods: it always converged while the other methods quickly stopped working when $\alpha$ becomes large. Hence it is observed that the phenomenon of power flow fractals is not exhibited by the proposed method unlike NR based methods \cite{thorp97fractals,Klump2000}. This hints that the fixed point method may avoid being trapped in local optima that can impact descent algorithms since local optima are not fixed points by definition.

%% file: conclusion.tex
A new fixed-point formulation of the power flow equation is developed in this paper. In contrast to existing fixed point formulations, it includes all possible cases of PV/PQ buses, mesh networks, resistive and inductive lines. Geometrically, our formulation treats the active and reactive power flow equations in rectangular voltage coordinates as circles and the power flow solutions as the intersection of these circles. Using a 3-tuple vector representation of circles, we derive simple, efficient and numerically stable formulas to find their intersection points. Based on the fixed point equations, we develop a natural iterative fixed point algorithm to solve the power flow problem. Numerical studies on the standard IEEE benchmarks show that our algorithm is able to converge when other state-of-the-art algorithms diverges or becomes unstable. We also show that the performance of proposed algorithm is comparable to other Jacobian based methods for large test systems. In addition, we show that our algorithm is robust to the initial starting point, able to converge for a wide range of starting conditions while other algorithms diverged.
% he numerical section shows that the proposed power flow method can converge for both base case and heavy loading. It is also observed that the proposed method can solve the problem when Jacobian becomes ill-conditioned. Finally, the sensitivity of the proposed method for random initialization is found to be insensitive for a PQ bus network. 

%% file: app.tex
%% !TEX root=main.tex
\subsection{Handling PQ and PV buses} \label{subsec:handle_PQ_PV_bus}
Here we present the fixed-point algorithm for a system with both PQ and PV buses. In the case of PV buses, Section~\ref{subsec:PV_bus} discusses the replacement of the reactive power balance equation by a condition on the voltage magnitude \eqref{eqn:voltage}. The voltage circle is centered at origin ($\mathit{c}_{v}$) with fixed radius of $\mathit{V}_{ref}$ ($\mathit{r}_{v}$). Thus the voltage solution for a PV bus can be calculated similarly to PQ bus by the intersection of two circles, real power \eqref{math_p} and specified voltage magnitude circles \eqref{eqn:voltage} at bus $d$.
\subsection{Bus Type Switching for PV Buses}
\subsubsection{PV to PQ switching} When there are no common points between the real power and specified voltage magnitude circles, the reactive power at bus $d$ is calculated to check for the violation of reactive power limits. In such a scenario, the PV bus is converted to a PQ bus by fixing its reactive power with the violated limit and it is then solved as a PQ bus. During the iterative process, this PQ bus is converted back to a PV bus as discussed below.
\subsubsection{PQ to PV switching} The bus that is converted to PQ has its real and reactive power fixed while the voltage phase angle and magnitude are free to change. However, this bus should be reverted to PV bus during the power flow iterative process when it is feasible since the problem still didn't converge. Let the voltage solution at this converted PQ bus be $v$. The violated upper and lower limit reactive powers be represented by $Q_{max}$ and $Q_{min}$ respectively. A converted PQ bus (due to PV to PQ switching) with $Q=Q_{max}$  or $Q=Q_{min}$ is referred as $PQ_{max}$ bus or $PQ_{min}$ bus respectively. For $PQ_{max}$ bus, when $v > \mathit{V}_{ref}$ it indicates that the reactive power at this bus is no longer needed to be fixed at $Q_{max}$ and it can be reverted back to a PV bus. Similarly for $PQ_{min}$ bus as well but the conditions to check is reversed to $v < \mathit{V}_{ref}$.

% \textcolor{blue}{Dear Dr.Zhang, I just want to let you know that I didn't code in the above approach discussed. Please critically read the above paragraph once for any mistakes, I thought a lot about the logic written above to make sure there are no mistakes if we use V circle.}
%
% \textcolor{blue}{I calculated the Q at a PV bus using neighboring bus voltages and calculated the voltage at the PV bus. Then I checked for QLimit violation. If not violated, I adjusted the voltage solution magnitude to be same as $V_{ref}$ at PV bus. If violated, I simply made it PQ bus and fixed Q to be Qmax or Qmin and solved it.}

% ref statement:
% Those buses called voltage regulated buses. ... When its reactive power injection reached its up- per or lower limit, the type of this bus becomes PQ, which means that the real and reactive power injections are fixed while the voltage phase angle and magnitude are free.

\begin{algorithm}[H]
\caption{RFPPF method for PQ and PV buses.}\label{alg:fixedpoint2}
\begin{algorithmic}[1]
%  \Procedure{myalgo}{}
\Statex \textbf{Input~~:} $\mathit{P}_i$ for bus $i=2,\cdots,n$,
\Statex ~~~~~~~~~~~$\mathit{Q}_i$ $\forall$ PQ buses,
\Statex ~~~~~~~~~~~$\mathit{V}_{ref}$ $\forall$ PV buses,
\Statex ~~~~~~~~~~~$\mathit{Q}_{max}$ and $\mathit{Q}_{min}$ for PV buses,
\Statex ~~~~~~~~~~~~~Bus type information $B\{i\}$ for bus $i=2,\cdots,n$,
\Statex ~~~~~~~~~~~Tolerance $\delta$ for the stopping criterion.
\Statex \textbf{Output:} $\mathit{v}_i$ for bus $i=2,\cdots,n$.
%\vspace{-1em}
\State Initialize voltages at all buses, $\mathit{v}_i$ for $i=2,\cdots,n$;
\State  Let the neighboring bus index be $k$;
\Statex where $k$ $\in$ $\mathbb{N}_m$.
\Comment{$\mathbb{N}_m$}{: Neighboring buses of bus $m$.}
\State Calculate the power mismatch ($\Delta\mathbb{S}$);

\While{($\Delta\mathbb{S}) > \delta$}
\Comment{}{Convergence criteria.}
\For{$m = 2$, $\cdots$, $n$}
\If{$B\{m\}$ == PQ bus}
\LState Calculate $(\mathit{c}_{p},\mathit{r}_{p})$ at bus $m$, $\forall$ $k$ $\in$ $\mathbb{N}_m$;
\LState Calculate $(\mathit{c}_{q},\mathit{r}_{q})$  at bus $m$, $\forall$ $k$ $\in$ $\mathbb{N}_m$;
\LState Calculate the voltage $\mathit{v}_m$ for bus $m$;
\LState $\mathit{v}_m = \mathit{v}_m$;
\Comment{}{Update bus $m$ voltage.}
\If{$B\{m\} == PQ_{max}$ and $\mathit{v}_m > \mathit{V}_{ref}$}
\LState Revert to PV bus and update $B\{m\}$;
\EndIf
\If{$B\{m\} == PQ_{min}$ and $\mathit{v}_m < \mathit{V}_{ref}$}
\LState Revert to PV bus and update $B\{m\}$;
\EndIf
\Else \Comment{}{PV bus.}
\LState Calculate $(\mathit{c}_{p},\mathit{r}_{p})$ at bus $m$, $\forall$ $k$ $\in$ $\mathbb{N}_m$;
\LState Calculate $(\mathit{c}_{v},\mathit{r}_{v})$ at bus $m$, $\forall$ $k$ $\in$ $\mathbb{N}_m$;
\If{Circles $(\mathit{c}_{p},\mathit{r}_{p})$ and $(\mathit{c}_{v},\mathit{r}_{v})$ intersect}
\LState Calculate the voltage $\mathit{v}_m$ for bus $m$;
\LState $\angle\mathit{v}_m= \angle\mathit{v}_m$;
\Comment{}{Update voltage angle.}
\Else
\LState Convert bus $m$ to PQ bus \& update $B\{m\}$;
\LState Calculate $(\mathit{c}_{p},\mathit{r}_{p})$  at bus $m$;
\LState Calculate $(\mathit{c}_{q},\mathit{r}_{q})$  at bus $m$;
\LState Calculate the voltage $\mathit{v}_m$ for bus $m$;
\LState $\mathit{v}_m = \mathit{v}_m$;
\Comment{}{Update bus $m$ voltage.}
\EndIf

\EndIf

\EndFor
\LState Calculate ($\Delta\mathbb{S}$);
\EndWhile \\
\Return $\mathit{v}_m$ $\forall$ buses $m=2,\cdots,n$;
%  \EndProcedure
\end{algorithmic}
\end{algorithm}